\documentstyle[12pt,amsbsy,amscd,amsfonts,amssymb,amstex,amstext,epsf]{amsart}
\theoremstyle{plain}
   \newtheorem{thm}{Theorem}[section]
   \newtheorem{lem}[thm]{Lemma}
   
   \newtheorem{defn-prop}[thm]{Definition-Proposition}
\theoremstyle{definition}
   \newtheorem{defn}[thm]{Definition}
   
   \newtheorem{prob}[thm]{Problem}
   \newtheorem{exmp}[thm]{Example}
\theoremstyle{remark}
   \newtheorem{rem}[thm]{Remark}
   \newtheorem{claim}[thm]{Claim}
   \newtheorem{rec}[thm]{Recall}

\begin{document}
\title{On $\delta_m$ Constant Locus of Versal Deformations of \\
        Nondegenerate Hypersurface Simple K3 Singularities}
\author{Masako Furuya}
\email{furuya@@math.titech.ac.jp}
\maketitle
\section*{Introduction}

   Simple $K3$ singularities are regarded as natural generalizations in 
three-dimensional case of simple elliptic singularities. The notion of a 
simple $K3$ singularity was defined by S. Ishii and K. Watanabe [IW] as a 
three-dimensional Gorenstein purely elliptic singularity of (0,2)-type, 
whereas a simple elliptic singularity is a two-dimensional purely elliptic 
singularity of (0,1)-type. It is also pointed out in [IW] that a simple $K3$ 
singularity is characterized as a quasi-Gorenstein singularity such that the 
exceptional set of any minimal resolution is a normal $K3$ surface. 
Let $ f \in {\bold C}[x,y,z,w] $ be a polynomial which is nondegenerate with 
respect to its Newton boundary $\Gamma(f)$ in the sense of [V1], and whose zero
locus $ X=\{f=0\} $ in ${\bold C}^4$ has an isolated singularity at the origin 
$ 0 \in {\bold C}^4 $. Then the condition for $(X,0)$ to be a simple $K3$ 
singularity is given by a property of the Newton boundary $\Gamma(f)$ of $f$ 
(cf. Theorem 1.6). 

   Hypersurface simple $K3$ singularities defined by nondegenerate 
quasi-homogeneous polynomials are classified into ninety five classes in term 
of the weight of the polynomial by Yonemura [Yo]. We consider versal 
deformations of them. It has been conjectured that the stratum $\mu$ =const of 
the versal deformation of any nondegenerate hypersurface simple $K3$ 
singularity is equivalent to the $\delta_m$ constant locus by Ishii. 
It holds true for the case deformations are also nondegenerate by 1.7 (1) [W]. 
On the other hand, it follows from 2.2 ([R1], [R2]) that the $\delta_m$ 
constant locus includes the $\mu$ constant locus generally. We show the 
conjecture holds true in general for No.10-14, 46-51 and 83 in the table of 
[Yo]. 

   I would like to express my sincere gratitude to Professor Shihoko Ishii for
telling me about the conjecture and giving valuable advice. I also express my 
gratitude to \\
Professors Masataka Tomari and Kei-ichi Watanabe for their useful comments \\
concerning Theorem 2.4, and to Professor Takao Fujita for his helpful remark 
on \\ Section 1 and 2. I also thank Mr. Hironobu Ishihara who pointed out 
grammatical mistakes throughout this paper. 
\section{Preliminary}

   In this section, we recall some definitions and facts from [I1], [IW], [W] 
and [Yo].

   First we define the plurigenera $ \delta_m, \; m \in \bold N $, for normal 
isolated singularities and define purely elliptic singularities. Let $(X,x)$ 
be a normal isolated singularity in an $n$-dimensional analytic space $X$, and 
$ \pi \; : \; (\tilde{X},E) \longrightarrow (X,x) $ a good resolution. In the 
following, we assume that $X$ is a sufficiently small Stein neighbourhood of 
$x$. 
\begin{defn}[Watanabe {\rm [W]-Def. 1.2}]
Let $(X,x)$ be a normal isolated singularity. For any positive integer $m$, 
$$ \delta_m(X,x):=\dim_{\bold C}\Gamma(X-\{x\},{\cal O}(mK))/L^{2/m}(X-\{x\})
, $$
where $K$ is the canonical line bundle on $ X-\{x\} $, and $ L^{2/m}(X-\{x\}) $
is the set of all $L^{2/m}$-integrable (at $x$) holomorphic $m$-ple $n$-forms 
on $ X-\{x\} $. 
\end{defn}

   Then $ \delta_m $ is finite and does not depend on the choice of a Stein 
neighbourhood $X$. 
\begin{defn}[Watanabe {\rm [W]-Def. 3.1}]
A singularity $(X,x)$ is said to be purely elliptic if $ \delta_m(X,x)=1 $ for 
every $ m \in \bold N $.
\end{defn}

   In the following, we assume that $(X,x)$ is quasi-Gorenstein, i.e., there 
exists a non-vanishing holomorphic 3-form on $ X-\{x\} $. Let $ E=\bigcup E_i $
be the decomposition of the exceptional set $E$ into irreducible components, 
and write 
$ K_{\tilde{X}} = \pi^{\ast}K_X+\sum_{i \in I}m_i E_i-\sum_{j \in J}m_j E_j $ 
with $ m_i \geq 0, \; m_j > 0 $. Ishii [I1] defined the essential part of the 
exceptional set $E$ as $ E_J = \sum_{j \in J}m_j E_j $, and showed that if 
$(X,x)$ is purely elliptic, then $ m_j=1 $ for all $ j \in J $. 
\begin{defn}[Ishii {\rm [I1]-Def. 4.1}]
A quasi-Gorenstein purely elliptic singularity $(X,x)$ is of $(0,i)$-type if 
$ H^{n-1}(E_J, {\cal O}_{E_J}) $ consists of the $(0,i)$-Hodge component 
$ H_{n-1}^{0,i}(E_J) $, where 
$$ {\bold C} \cong H^{n-1}(E_J, {\cal O}_{E_J}) = Gr_F^0 H^{n-1}(E_J) 
   = \bigoplus_{i=0}^{n-1} H_{n-1}^{0,i}(E_J). $$
\end{defn}
\begin{defn-prop}[Ishii-Watanabe {\rm [IW]-Def. 4}]
A three-dimensional \\ singularity $(X,x)$ is a simple $K3$ singularity 
if the following two equivalent conditions are satisfied\rom:

  \rom{(1)} \ \ $(X,x)$ is Gorenstein purely elliptic of (0,2)-type.

  \rom{(2)} \ \ The exceptional divisor $E$ is a normal $K3$ surface for any 
minimal resolution $ \pi : (\tilde{X},E) \longrightarrow (X,x) $. 
\end{defn-prop}
\begin{rec}
A minimal resolution $ \pi : (\tilde{X},E) \longrightarrow (X,x) $ is a proper
morphism with $ \tilde{X}-E \cong X-\{x\} $, where $\tilde{X}$ has only 
terminal $ \bold Q $-factorial singularities and $K_{\tilde{X}}$ is numerically
effective with respect to $\pi$.
\end{rec}

   Next we consider the case where $(X,x)$ is a hypersurface singularity
defined by a nondegenerate polynomial $f= \sum a_{\nu} z^{\nu} \in 
{\bold C}[z_0,z_1,\dotsb,z_n]$, and $x=0 \in {\bold C}^{n+1}$. 
\begin{rec}
The Newton boundary $\Gamma(f)$ of $f$ is the union of the compact faces of 
$\Gamma_+(f)$, where $\Gamma_+(f)$ is the convex hull of $ \bigcup_{a_{\nu} \ne
0}(\nu+{\bold R}_{\geq 0}^{n+1}) $ in ${\bold R}^{n+1}$. For any face $\Delta$
of $\Gamma_+(f)$, set $ f_{\Delta}:=\sum_{\nu \in \Delta}a_{\nu}z^{\nu} $. 
We say $f$ to be nondegenerate, if 
$$ \displaystyle \frac{\partial f_{\Delta}}{\partial z_0} = 
                 \frac{\partial f_{\Delta}}{\partial z_1} = \dotsb = 
                 \frac{\partial f_{\Delta}}{\partial z_n} = 0 $$
has no solution in $ ({\bold C}-\{0\})^{n+1} $ for any face $\Delta$. 
\end{rec}

   When $f$ is nondegenerate, the condition for $(X,x)$ to be a purely 
elliptic singularity is given as follows: 
\begin{thm}[Watanabe {\rm [W]-Prop. 2.9, Cor. 3.14}]
Let $f$ be a nondegenerate \\ polynomial and suppose $X=\{f=0\}$ has an 
isolated singularity at $x=0 \in {\bold C}^{n+1}$.

  \rom{(1)} $(X,x)$ is purely elliptic if and only if $(1,1,\dotsb,1) \in 
\Gamma(f)$.

  \rom{(2)} Let n=3 and let $\Delta_0$ be the face of $\Gamma (f)$ containing 
$(1,1,1,1)$ in the relative interior of $\Delta_0$. Then $(X,x)$ is a simple 
$K3$ singularity if and only if $\dim_{\bold R} \Delta_0=3$. 
\end{thm}

   Thus if $f$ is nondegenerate and defines a simple $K3$ singularity, then 
$f_{\Delta_0}:=\sum_{\nu \in \Delta_0} a_{\nu} z^{\nu}$ is a quasi-homogeneous
polynomial of a uniquely determined weight $\alpha$ called the weight of $f$ 
and denoted $\alpha(f)$. Namely, $\alpha =(\alpha_1,\dotsb,\alpha_4) \in {\bold
Q_{>0}}^4$ and $\deg_{\alpha}(\nu):=\sum_{i=1}^4 \alpha_i \nu_i =1$ for any 
$\nu \in \Delta_0$. In particular, $\sum_{i=1}^4 \alpha_i=1$, since $(1,1,1,1)$
is always contained in $\Delta_0$. 
\begin{thm}[Yonemura {\rm [Yo]-Prop. 2.1}]
The cardinality of \\ 
$ \{ \alpha (f) \; ; f \text{ is nondegenerate and
                     defines a simple } K3 \text{ singularity}, \; 
                     \alpha_1 \geq \alpha_2 \geq \alpha_3 \geq \alpha_4 \} $ \\
is 95. 
\end{thm}

  In Table 2.2 of [Yo] can be found the complete list of weights 
$\alpha=\alpha(f)$ and examples of $f=\sum \alpha_{\nu} z^{\nu}$ such that 
$f$ is quasi-homogeneous and that $\{f=0\} \subset {\bold C}^4$ has a simple 
$K3$ singularity at the origin $0 \in {\bold C}^4$.

   We describe a weight $\alpha=\alpha(f)$ as $\alpha=(p_1/p,p_2/p,p_3/p,
p_4/p)$, where $p, \; p_i$ are \\ positive integers with 
$ {\bold {gcd}}(p_1,p_2,p_3,p_4)=1.$ \\

   Next we consider the versal deformation of simple $K3$ singularity 
$(X=\{f=0\}, \; 0)$ such that $f$ is a nondegenerate quasi-homogeneous 
polynomial. 
\begin{defn}
A deformation of an isolated singularity $(X,x)$ is a flat family of 
singularities:
$$ {\frak X} =\{(X_t,x_t) \; ; \; t \in U \subset {\bold C}^N \} 
   \overset{flat}{\longrightarrow} U   $$
such that $(X_0,x_0) \cong (X,x)$ as germs of holomorphic functions. \\
Where $U$ is a sufficiently small open neighbourhood of 0 in $ {\bold C}^N $.
\end{defn}
\begin{defn}
A deformation $\frak X$ of an isolated singularity $(X,x)$ is versal if for 
any deformation ${\frak X'} \longrightarrow U'$ of $(X,x)$, the following is 
satisfied: 
\begin{equation*}
  \begin{CD}
    {\frak X}{\times}_U U' \overset{isom}{\cong} {\frak X}' @>>> U' \\
    @.
        @VV{\exists \; holomorphic}V \\
    {\hspace{2.1cm} \frak X} @>>> U \\
  \end{CD} 
\end{equation*}
\end{defn}
\begin{thm}[{\rm [KS], [Tj]}]
The versal deformation $\frak X$ of an isolated singularity \\ 
$(\{f=0\}, \; 0)$ is described by 
$$ {\frak X} =\{(\{f+\sum \lambda_i g_i=0\}, \; 0) \; ; 
                \; (\lambda_i) \in U \subset {\bold C}^N \}, $$
where the $g_i$ determine a $\bold C$-basis of the vector space 
${\bold C}\{z_0,\dotsb,z_n\}/(f, f_{z_0}, \dotsb, f_{z_n})$. 
\end{thm}

   The purpose of this paper is to show that the following conjecture holds 
true for No.10-14, 46-51 and 83 in Table 2.2 of [Yo].
\begin{prob}[Ishii]
Let $(X=\{f=0\}, \; 0)$ be a hypersurface simple $K3$ singularity 
defined by a nondegenerate quasi-homogeneous polynomial $f$, and let 
$\frak X$ be the versal deformation of $(X,0)$. Then, \\
\hspace*{1.5cm} $ \{ \lambda \in U \subset {\bold C}^N \; ; \; 
                     \mu(X,0)=\mu(X_{\lambda} ,0), \; 
                     (X_{\lambda} ,0) \in {\frak X} \} $ \\
\hspace*{1cm} $ = \{ \lambda \in U \subset {\bold C}^N \; ; \;
                  1=\delta_m (X,0)=\delta_m (X_{\lambda} ,0) \; 
                  \text{for all } m \geq 1, \; 
                  (X_{\lambda} ,0) \in {\frak X}  \}. \qquad \square $ 
\end{prob}

   Since $(X,0)$ and $(X_{\lambda} ,0)$ are hypersurface isolated 
singularities, they are normal Gorenstein, and so 
$$ P_g=0 \; \Longleftrightarrow \; \delta_m=0 \text{\; for all \ } m \geq 1.$$

   On the other hand, $(X_{\lambda} ,0)$ is a deformation of a purely elliptic
singularity $(X,0)$, so it is either rational or purely elliptic. 
Therefore, this problem is equivalent to:
$$   \mu(X,0)=\mu(X_{\lambda} ,0) \quad \Longleftrightarrow \quad 
     1=P_g(X,0)=P_g(X_{\lambda} ,0). $$

   Furthermore, since $\mu$ and $P_g$ are upper semi-continuous in respect of 
deformation \\([M], [Te], [E]-Thm. 1, [Ya]-Thm. 2.6), it is equivalent to:
$$      \mu (X,0) > \mu (X_{\lambda} ,0) \quad \Longleftrightarrow \quad
        1 = P_g (X,0) > P_g (X_{\lambda} ,0) = 0.    $$
\section{Reduction of the problem}

   Considering some facts, we can reduce the problem posed in section 1 to one
about the weight of the defining polynomial of a hypersurface singularity. 
\begin{thm}[Varchenko {\rm [V2]-Thm. 2}]
Let $f \in {\bold C}[z_0,\dotsb,z_n]$ be quasi-homogeneous of weight 
$ \alpha $ with $ \alpha_0+\dotsb+\alpha_n=1 $, and $\{f=0\}$ has an isolated 
singularity at $0 \in {\bold C}^{n+1}$. Let
$$     \mu := \mu (f,0)<+\infty, \qquad
        f_{\lambda} := f+\sum_{i=1}^{\mu} \lambda_i g_i, $$
where $g_i \in {\bold C}\{z_0,\dotsb,z_n\}/(f_{z_0},\dotsb,f_{z_n})$ are 
generators of the Jacobi ring, which are monomials. \rom(Since $f$ is 
quasi-homogeneous, $g_i$ can be taken as monomials.\rom) Then, \\
\hspace*{1.5cm} $ \{ \lambda =(\lambda_1,\dotsb,\lambda_{\mu} ) \in U \subset 
                   {\bold C}^{\mu} \; ; \; \mu (f,0)=\mu (f_{\lambda},0)\} $ \\
\hspace*{1cm} $ = \{ \lambda =(\lambda_1,\dotsb,\lambda_{\mu} ) \in U \subset 
                     {\bold C}^{\mu} \; ; \; \lambda_i=0 \text{\ for all i 
                     satisfying that } \deg_{\alpha} (g_i)<1 \}, $ \\
where $ \deg_{\alpha}(z_0^{i_0} \cdot \dotsb \cdot z_n^{i_n})
      :=\alpha_0 i_0+ \dotsb +\alpha_n i_n. \qquad \square $ 
\end{thm}

   From 1.7 (1) and 2.1, it follows that
$$ \mu(f,0)= \mu(f_{\lambda},0) \quad \Longleftrightarrow \quad
   P_g(f,0)=P_g(f_{\lambda},0) $$
holds true if $f_{\lambda}$ is also nondegenerate. Though $f_{\lambda}$ is not 
always nondegenerate, the \\ following theorem is useful as well. 
\begin{thm}[Reid {\rm [R1]-Thm. 4.1, [R2]-Thm. 4.6}]
Let $(X=\{f=0\}, \; 0) \subset ({\bold C}^{n+1},0)$ be a hypersurface 
singularity and let $\alpha=(\alpha_0,\dotsb,\alpha_n)=(p_0/p,\dotsb,p_n/p) 
\in {{\bold Q}_{>0}}^{n+1}$ such that $(p_0,\dotsb,p_n) \in {\bold N}^{n+1}$ 
is a primitive vector. Then, \\
\hspace*{2cm} $ (X,0): \; \text{canonical} \quad \Longrightarrow \quad
                \deg_{\alpha}(z_0\cdot\dotsb\cdot z_n)>\deg_{\alpha}(f), $ \\
where $ \deg_{\alpha}(f):=\min\{\deg_{\alpha}(z^{\nu}) \; ; 
                               \; z^{\nu} \in f \}. \qquad \square $ 
\end{thm}

   A hypersurface singularity is canonical if and only if it is rational, 
and so it follows that
$$  \mu(f,0)=\mu(f_{\lambda},0) \Longrightarrow P_g(f,0)=P_g(f_{\lambda},0) $$
holds always true from 2.2. Hence we should show the converse proposition. 
\begin{rem}
Let $(X,x)$ be a n-dimensional normal Gorenstein singularity and 
$$  (\Tilde{\Tilde X},\tilde{E}) \overset{\pi'}{\longrightarrow} (\tilde{X},E)
     \overset{\pi}{\longrightarrow} (X,x) , $$
where $\tilde{X}$ has at most rational singularities and 
$\pi'':=\pi \circ \pi'$ is a resolution. 
Let $ E=\bigcup_{i \in I}E_i $ be the decomposition of the exceptional set $E$
into irreducible components, and 
write $ K_{\tilde{X}}=\pi^{\ast} K_X+\sum_{i \in I} m_i E_i $. 
Then, it follows that 
$ \pi'_{\ast}({\cal O}(K_{\Tilde{\Tilde X}}))={\cal O}(K_{\tilde{X}}) $ since 
$\tilde{X}$ has at most rational singularities ([KKMS]-p. 50), and so
{\allowdisplaybreaks
\begin{align*}
  P_g(X,x) & :=\dim_{\bold C}
               (R^{n-1}{\pi''}_{\ast} {\cal O}_{\Tilde{\Tilde X}})_x  \\
           &  =\dim_{\bold C}\Gamma(X-\{x\},{\cal O} (K_X))/L^2(X-\{x\}) \\
           &  =\dim_{\bold C}\Gamma(\Tilde{\Tilde X}-\tilde{E}, 
                                   {\cal O} (K_{\Tilde{\Tilde X}}))/
                             \Gamma(\Tilde{\Tilde X}, 
                                   {\cal O} (K_{\Tilde{\Tilde X}})) \\
           &  =\dim_{\bold C}
               \Gamma(\Tilde{\Tilde X}, {\cal O} ({\pi''}^{\ast} K_X))/
               \Gamma(\Tilde{\Tilde X}, {\cal O} (K_{\Tilde{\Tilde X}})) \\
           &  =\dim_{\bold C}\Gamma(\tilde{X}, {\cal O} (\pi^{\ast} K_X))/
              \Gamma(\tilde{X}, \pi'_{\ast}({\cal O}(K_{\Tilde{\Tilde X}}))) \\
           &  =\dim_{\bold C}\Gamma(\tilde{X}, {\cal O} (\pi^{\ast} K_X))/
                             \Gamma(\tilde{X}, {\cal O} (K_{\tilde{X}})). 
\end{align*}}
Therefore,
$$  \exists \; i \in I \; ; \; m_i<0 \; \Longleftrightarrow \;
    \Gamma(\tilde{X},{\cal O} (\pi^{\ast}K_X)) \supsetneq 
    \Gamma(\tilde{X},{\cal O} (K_{\tilde{X}})) \;
    \Longleftrightarrow \; P_g(X,x) > 0, $$
namely,
$$ \forall \; i \in I , \; m_i \geq 0 \; \Longleftrightarrow \; P_g(X,x)=0. $$
\end{rem}
\begin{thm}[Tomari-Watanabe {\rm [TW]-Thm. 5.6}]
Let $(X=\{f=0\}, \; x)$ be a \\ n-dimensional hypersurface isolated 
singularity and \\
\hspace*{1cm} $ f=f_0+f_1+f_2+\dotsb, $ \\
\hspace*{2cm} $ f_i : $ quasi-homogeneous polynomial of weight
                \; $ \alpha=(\alpha_0,\dotsb,\alpha_n) \; ; $ \\
\hspace*{2cm} $ 1=\deg_{\alpha}(f_0)<\deg_{\alpha}(f_1)<\deg_{\alpha}(f_2)
                 <\dotsb, $ \\
and \\
\begin{figure}[h]
\setlength{\unitlength}{1mm}
\begin{picture}(117,41)(-35,0)
  \put(0,0){\makebox(45,7){$ E={\pi}^{-1}(x) $}}
  \put(45,0){\makebox(25,7){$ \longrightarrow $}}
  \put(70,0){\makebox(12,7){$ x $.}}
  \put(0,7){\makebox(45,9){$ \bigcup $}}
  \put(70,7){\makebox(12,9){$ \cup $}}
  \put(76,10.5){\line(0,1){2.5}}
  \put(0,16){\makebox(45,7)
            {$ \tilde{X}=\overline{{\Pi}^{-1}(X)-{\Pi}^{-1}(x)} $}}
  \put(45,16){\makebox(25,9)
             {$ @>{\pi=\Pi {\vert}_{\tilde{X}}}>> $}}
  \put(70,16){\makebox(12,7){$ X $}}
  \put(0,23){\makebox(45,9){$ \bigcup $}}
  \put(70,23){\makebox(12,9){$ \bigcup $}}
  \put(0,32){\makebox(45,7){$ V $}}
  \put(45,32){\makebox(25,9)
             {$ @>{\Pi \; : \; \alpha-blow-up}>> $}}
  \put(70,32){\makebox(12,7){$ {\bold C}^{n+1} $}}
\end{picture}
\end{figure} \\
Assume that $f_0$ is irreducible and both $X-\{x\}$ and $\{f_0=0\}-\{x\}$ have
at most rational singularities around $x$. Then $\tilde{X}$ has at most 
rational singularities. \qquad $\square$ 
\end{thm}

  Thus we expect a partial resolution $\pi$ in 2.3 is given by a weighted 
blow-up. 
\begin{rem}[Ishii {\rm [I2]-Prop. 1.3, 1.6}]
Under the notation in 2.4, $ \Pi^{\ast} X \subset V $ and $ K_{{\bold C}^4} $ 
are principal divisors, hence,
{\allowdisplaybreaks
\begin{align*}
 \Pi^{\ast} X & = \tilde{X} + p F, \\
 K_V & = \Pi^{\ast} K_{{\bold C}^4} + (p_1+p_2+p_3+p_4-1)F, \\
 K_X & = (K_{{\bold C}^4}+X)|_X, \\
 \text{thus,} \qquad
 K_{\tilde{X}} & = (K_V+\tilde{X})|_{\tilde{X}} \\
               & = (\Pi^{\ast}(K_{{\bold C}^4}+X) + (p_1+p_2+p_3+p_4-1-p) F)|_
                   {\tilde{X}} \\
               & = \pi^{\ast}K_X+(p_1+p_2+p_3+p_4-1-p)\sum_i k_i E_i,
\end{align*}}
where \ $ \alpha = (\alpha_1, \; \alpha_2, \; \alpha_3, \; \alpha_4)
                 = (p_1/p, \; p_2/p, \; p_3/p, \; p_4/p) $ with 
        $ {\bold {gcd}}(p_1, \; p_2, \; p_3, \; p_4)=1 $,\\
\hspace*{1.3cm} $ F={\Pi}^{-1}(0), \quad 
                  F|_{\tilde{X}}=E=\sum_i k_i E_i, \quad k_i>0. $ \\
Thus, \ \ $ p<p_1+\dotsb+p_4 $ if and only if
          $ R^2 \pi_{\ast} {\cal O}_{\tilde{X}}=0 $.
\end{rem}
\begin{prob}
Let $f$ be a nondegenerate quasi-homogeneous polynomial which \\ defines a 
simple $K3$ singularity at $0$, and $(X_{\lambda}=\{f_{\lambda}=0\}, \; 0)$
a versal deformation of \\ $(X=\{f=0\}, \; 0)$ such that $ \mu(X_{\lambda},0) 
< \mu(X,0) $. Then find a weight $\alpha'=(\alpha'_1,\dotsb,\alpha'_4)
=(p'_1/p',\dotsb,p'_4/p')$ of $f_{\lambda}$ with 
$ 1<\alpha'_1+\dotsb+\alpha'_4 $ such that: \\
\hspace*{1cm} $ f_{\lambda}=f_0+f_1+f_2+\dotsb, $ \\ 
\hspace*{2cm} $ f_i $ : quasi-homogeneous polynomial of weight $ \alpha', $ \\
\hspace*{2cm} $ 1=\deg_{\alpha'}(f_0)<\deg_{\alpha'}(f_1)<\deg_{\alpha'}(f_2)<
                \dotsb, $ \\
\hspace*{2cm} $ f_0 $ is irreducible, \\
\hspace*{1cm} and $ \{f_0=0\}-\{0\} $ has at most rational singularities around
              $ 0 $. 
\end{prob}

   If there exists such weight $\alpha'$ then $P_g(X_{\lambda},0)=0$ by 2.4 
and 2.5. 

   We show there exists such weight $\alpha'$ as in 2.6 for 
$ f=x^2+y^3+\dotsb, $ of No.10-14, 46-51 and 83 in Table 2.2 of [Yo] to obtain 
our main result as follows: 
\begin{thm}
Let $(X=\{f=0\}, \; 0)$ be a hypersurface simple $K3$ singularity 
defined by a nondegenerate quasi-homogeneous polynomial $ f=x^2+y^3+\dotsb, $ 
which is one of No.10-14, 46-51 and 83 in Table 2.2 of [Yo]\rom: \\
\hspace*{1.5cm} $($No. 10$)$ \ \ $ f = x^2 + y^3 + z^{12} + w^{12} $,\\
\hspace*{1.5cm} $($No. 11$)$ \ \ $ f = x^2 + y^3 + z^{10} + w^{15} $,\\
\hspace*{1.5cm} $($No. 12$)$ \ \ $ f = x^2 + y^3 + z^9 + w^{18} $,\\
\hspace*{1.5cm} $($No. 13$)$ \ \ $ f = x^2 + y^3 + z^8 + w^{24} $,\\
\hspace*{1.5cm} $($No. 14$)$ \ \ $ f = x^2 + y^3 + z^7 + w^{42} $,\\
\hspace*{1.5cm} $($No. 46$)$ \ \ $ f = x^2 + y^3 + z^{11} + z w^{12} $,\\
\hspace*{1.5cm} $($No. 47$)$ \ \ $ f = x^2 + y^3 + y z^7 + z^9 w^2+w^{14}$,\\
\hspace*{1.5cm} $($No. 48$)$ \ \ $ f = x^2 + y^3 + z^9 w + w^{16} $,\\
\hspace*{1.5cm} $($No. 49$)$ \ \ $ f = x^2 + y^3 + z^8 w + w^{21} $,\\
\hspace*{1.5cm} $($No. 50$)$ \ \ $ f = x^2 + y^3 + y z^5 + z^7 w^2+w^{30}$,\\
\hspace*{1.5cm} $($No. 51$)$ \ \ $ f = x^2 + y^3 + z^7 w + w^{36} $,\\
\hspace*{1.5cm} $($No. 83$)$ \ \ $ f = x^2 + y^3 + y w^9+z^{10} w+z^2w^{11}$,\\
and let $\frak X$ be the versal deformation of $(X,0)$. Then, \\
\hspace*{1.5cm} $ \{ \lambda \in U \subset {\bold C}^N \; ; \; 
                     \mu(X,0)=\mu(X_{\lambda} ,0), \; 
                     (X_{\lambda} ,0) \in {\frak X} \} $ \\
\hspace*{1cm} $ = \{ \lambda \in U \subset {\bold C}^N \; ; \;
                  1=\delta_m (X,0)=\delta_m (X_{\lambda} ,0) \; 
                  \text{for all } m \geq 1, \; 
                  (X_{\lambda} ,0) \in {\frak X}  \}. $ 
\end{thm}
\section{Proof of Theorem 2.7}

   We prepare the following Lemma prior to the proof of Theorem 2.7. 
\begin{lem}
Let $f$ be a nondegenerate quasi-homogeneous polynomial listed by \\ Yonemura 
[Yo] which defines a simple $K3$ singularity, and $f_{\lambda}$ a versal 
deformation of $f$ such that $ \mu(f_{\lambda},0) < \mu(f,0) $. 
Under this original local coordinate system, if 
$ \alpha'=(\alpha'_1,\dotsb,\alpha'_4) \in {{\bold Q}_{>0}}^4 $ satisfies 
$ \deg_{\alpha'}(f_{\lambda})=1 $ then $ 1<\alpha'_1+\dotsb+\alpha'_4. $ 
\end{lem}
\begin{pf}

   For any i=1,2,3,4, one of the following is satisfied:
 
   (a) \ $ p_i | p $,

   (b) \ $ p_i | (p-p_j) $ for some $ j \neq i $,\\
namely, 

   (a) \ $ z_i^I \in f \qquad (\exists I \geq 2) $,

   (b) \ $ z_i^I z_j \in f \qquad (\exists I \geq 2) $ for some $ j \neq i $.\\
By assumptions, $ p \alpha_1+q \alpha_2+r \alpha_3+s \alpha_4<1 $ 
                for some $ {z_1}^p {z_2}^q {z_3}^r {z_4}^s \in f_{\lambda} $,\\
\hspace*{2.2cm} and $ p \alpha'_1+q \alpha'_2+r \alpha'_3+s \alpha_4 \geq 1 $ 
                for all $ {z_1}^p {z_2}^q {z_3}^r {z_4}^s \in f_{\lambda} $,\\
so there exists $ i \in \{1,2,3,4\} $ such that $ \alpha_i<\alpha'_i \qquad 
                \dotsb \circledast.$ \\ 
So the defining polynomials $f$ can be classified as below.\\ 
Case $\bold 1.$ \ \ $ z_i^I \in f \quad (\exists I \geq 2) $.\\
\hspace*{1cm} Then $ I \alpha_i=1 \leq I \alpha'_i $,  
              so we have $ \alpha_i \leq \alpha'_i. $ \\
\hspace*{0.2cm} ($\bold 1$-I) \ \ $ z_i^I, \; z_j^J \in f \quad
                                    (\exists I, \; \exists J \geq 2) $.\\
\hspace*{1.2cm} Since $ J \alpha_j=1 \leq J \alpha'_j $, 
                we have $ \alpha_j \leq \alpha'_j. $ \\
\hspace*{0.4cm} ($\bold 1$-I-i) \ \ $ z_i^I, \; z_j^J, \; z_k^K \in f \quad 
                            (\exists I, \; \exists J, \; \exists K \geq 2) $.\\
\hspace*{1.4cm} From $ K \alpha_k=1 \leq K \alpha'_k $, 
                we have $ \alpha_k \leq \alpha'_k. $ \\
\hspace*{0.6cm} ($\bold 1$-I-i-a) (No.1-14) \ \ 
                 $ z_i^I, \; z_j^J, \; z_k^K, \; z_l^L \in f \quad 
              (\exists I, \; \exists J, \; \exists K, \; \exists L \geq 2) $.\\
\hspace*{1.6cm} From $ L \alpha_l=1 \leq L \alpha'_l $, 
                we have $ \alpha_l \leq \alpha'_l $, 
                so $ \alpha_1+\dotsb+\alpha_4
                  <\alpha'_1+\dotsb+\alpha'_4 $ by $\circledast$. \\
\hspace*{0.6cm} ($\bold 1$-I-i-b) (No.15-51) \ \ 
                 $ z_i^I, \; z_j^J, \; z_k^K, \; z_k z_l^L \in f \quad 
              (\exists I, \; \exists J, \; \exists K, \; \exists L \geq 2) $.\\
\hspace*{1.6cm} From $ \alpha_k+L \alpha_l=1 \leq \alpha'_k+L \alpha'_l $, 
                we have $ \alpha_k+\alpha_l+(L-1)(\alpha_l-\alpha'_l) 
                          \leq \alpha'_k+\alpha'_l. $\\
\hspace*{1.6cm} If $ \alpha_l \leq \alpha'_l $ then $ \alpha_1+\dotsb+\alpha_4
                    <\alpha'_1+\dotsb+\alpha'_4 $ by $\circledast$.  \\
\hspace*{1.6cm} If $ \alpha_l>\alpha'_l $ then \\
\hspace*{1.6cm} $ \alpha _i+\alpha _j+\alpha _k+\alpha _l 
                < \alpha _i+\alpha _j+\alpha _k+\alpha _l
                  +(L-1)(\alpha_l-\alpha'_l) \leq 
                  \alpha'_i+\alpha'_j+\alpha'_k+\alpha'_l . $ \\
\hspace*{0.4cm} ($\bold 1$-I-ii) \ \ $ z_i^I, \; z_j^J, \; z_k^K z_l \in f 
                      \quad (\exists I, \; \exists J, \; \exists K \geq 2)$.\\
\hspace*{1.4cm} $ \alpha_k \leq \alpha'_k $ or $ \alpha_l \leq \alpha'_l $, \ 
                and $ \alpha_k+\alpha_l+(K-1)(\alpha_k-\alpha'_k) 
                      \leq \alpha'_k+\alpha'_l $.\\
\hspace*{0.6cm} ($\bold 1$-I-ii-a) (No.78) \ \ 
                 $ z_i^I, \; z_j^J, \; z_k^K z_l, \; z_k z_l^L \in f \quad 
              (\exists I, \; \exists J, \; \exists K, \; \exists L \geq 2) $.\\
\hspace*{1.6cm} $ \alpha_k+\alpha_l+(L-1)(\alpha_l-\alpha'_l) 
                  \leq \alpha'_k+\alpha'_l $.\\
\hspace*{1.6cm} If $ \alpha_k \leq \alpha'_k $ and $ \alpha_l \leq \alpha'_l $ 
                then $ \alpha_1+\dotsb+\alpha_4 <
                \alpha'_1+\dotsb+\alpha'_4 $ by $\circledast$.  \\
\hspace*{1.6cm} If $ \alpha_k \leq \alpha'_k $ and $ \alpha_l > \alpha'_l $ 
                then \\
\hspace*{1.6cm} $ \alpha _i+\alpha _j+\alpha _k+\alpha _l 
                < \alpha _i+\alpha _j+\alpha _k+\alpha _l
                  +(L-1)(\alpha_l-\alpha'_l) \leq 
                  \alpha'_i+\alpha'_j+\alpha'_k+\alpha'_l . $ \\
\hspace*{1.6cm} If $ \alpha_k > \alpha'_k $ and $ \alpha_l \leq \alpha'_l $ 
                then \\
\hspace*{1.6cm} $ \alpha _i+\alpha _j+\alpha _k+\alpha _l 
                < \alpha _i+\alpha _j+\alpha _k+\alpha _l
                  +(K-1)(\alpha_k-\alpha'_k) \leq 
                  \alpha'_i+\alpha'_j+\alpha'_k+\alpha'_l . $ \\
\hspace*{0.6cm} ($\bold 1$-I-ii-b) (No.52, 54-74, 76, 77, 79-83) \\
\hspace*{1.6cm} $ z_i^I, \; z_j^J, \; z_k^K z_l, \; z_j z_l^L \in f $ \quad
              $(\exists I, \; \exists J, \; \exists K, \; \exists L \geq 2)$.\\
\hspace*{1.6cm} $ \alpha_j+\alpha_l+(L-1)(\alpha_l-\alpha'_l) 
                  \leq \alpha'_j+\alpha'_l $. \ \ 
                The rest is similar to $\bold 1$-I-ii-a .\\
\hspace*{0.4cm} ($\bold 1$-I-iii) \ \ $ z_i^I, \; z_j^J, \; z_j z_k^K \in f 
                      \quad (\exists I, \; \exists J, \; \exists K \geq 2)$.\\
\hspace*{1.4cm} Similarly, $ \alpha_j+\alpha_k+(K-1)(\alpha_k-\alpha'_k) 
                             \leq \alpha'_j+\alpha'_k $.\\
\hspace*{0.6cm} ($\bold 1$-I-iii-a) (No.66, 67, 72, 75, 81, 82) \\
\hspace*{1.6cm} $ z_i^I, \; z_j^J, \; z_j z_k^K, \; z_j z_l^L \in f $ \quad
              $(\exists I, \; \exists J, \; \exists K, \; \exists L \geq 2)$.\\
\hspace*{1.6cm} Since there exists $ z_k^r z_l^s \in f $, we have  
                $ \alpha_k \leq \alpha'_k $ or $ \alpha_l \leq \alpha'_l $.\\
\hspace*{1.6cm} Moreover 
                $ \alpha_j+\alpha_l+(L-1)(\alpha_l-\alpha'_l) \leq 
                  \alpha'_j+\alpha'_l $ by $ z_j z_l^L \in f $.\\
\hspace*{1.6cm} So we have 
                $ \alpha_1+\dotsb+\alpha_4 < 
                \alpha'_1+\dotsb+\alpha'_4 $ similarly as in 
                $\bold 1$-I-ii-a.\\
\hspace*{0.6cm} ($\bold 1$-I-iii-b) (No.53, 57, 58, 62-64, 66, 67, 69-72) \\
\hspace*{1.6cm} $ z_i^I, \; z_j^J, \; z_j z_k^K, \; z_i z_l^L \in f $ \quad
              $(\exists I, \; \exists J, \; \exists K, \; \exists L \geq 2)$.\\
\hspace*{1.6cm} $ \alpha_i+\alpha_l+(L-1)(\alpha_l-\alpha'_l) \leq 
                  \alpha'_i+\alpha'_l $.\\
\hspace*{1.6cm} If $ \alpha_k \leq \alpha'_k $ or $ \alpha_l \leq \alpha'_l $ 
                then the rest is similar to $\bold 1$-I-ii-a. \\
\hspace*{1.6cm} If $ \alpha_k > \alpha'_k $ and $ \alpha_l > \alpha'_l $ then\\
\hspace*{1.6cm} $ \alpha_i+\alpha_j+\alpha_k+\alpha_l 
                < \alpha_i+\alpha_j+\alpha_k+\alpha_l+(K-1)(\alpha_k-\alpha'_k)
                  +(L-1)(\alpha_l-\alpha'_l) $ \\
\hspace*{4.8cm} $ \leq \alpha'_i+\alpha'_j+\alpha'_k+\alpha'_l .$ \\
\hspace*{0.2cm} ($\bold 1$-II) \ \ $ z_i^I, \; z_j^J z_k \in f \quad
                                    (\exists I, \; \exists J \geq 2) $.\\
\hspace*{1.2cm} $ \alpha_j \leq \alpha'_j $ or $ \alpha_k \leq \alpha'_k $, 
                \ and $ \alpha_j+\alpha_k+(J-1)(\alpha_j-\alpha'_j) \leq 
                        \alpha'_j+\alpha'_k $.\\
\hspace*{0.4cm} ($\bold 1$-II-i) \ \ $ z_i^I, \; z_j^J z_k, \; z_k^K z_l \in f 
                      \quad (\exists I, \; \exists J, \; \exists K \geq 2) $.\\
\hspace*{1.4cm} $ \alpha_k \leq \alpha'_k $ or $ \alpha_l \leq \alpha'_l $, 
                \ and $ \alpha_k+\alpha_l+(K-1)(\alpha_k-\alpha'_k) \leq 
                        \alpha'_k+\alpha'_l.$\\
\hspace*{0.6cm} ($\bold 1$-II-i-a) (No.88, 90-93) \ \ 
                $ z_i^I, \; z_j^J z_k, \; z_k^K z_l, \; z_k z_l^L \in f \quad
              (\exists I, \; \exists J, \; \exists K, \; \exists L \geq 2) $.\\
\hspace*{1.6cm} Since there exists $ z_j^q z_l^s \in f $, we have 
                $ \alpha_j \leq \alpha'_j $ or $ \alpha_l \leq \alpha'_l $.\\
\hspace*{1.6cm} Moreover $ \alpha_k+\alpha_l+(L-1)(\alpha_l-\alpha'_l) \leq 
                           \alpha'_k+\alpha'_l $ by $ z_k z_l^L \in f $.\\
\hspace*{1.6cm} If $ \alpha_j \leq \alpha'_j $, $ \alpha_k \leq \alpha'_k $ 
                and $ \alpha_l \leq \alpha'_l $, then the assertion holds true.
                \\
\hspace*{1.6cm} If $ \alpha_j \leq \alpha'_j $, $ \alpha_k \leq \alpha'_k $
                and $ \alpha_l > \alpha'_l $, then \\
\hspace*{1.6cm} $ \alpha _i+\alpha _j+\alpha _k+\alpha _l 
                < \alpha _i+\alpha _j+\alpha _k+\alpha _l
                  +(L-1)(\alpha_l-\alpha'_l) \leq 
                  \alpha'_i+\alpha'_j+\alpha'_k+\alpha'_l . $ \\
\hspace*{1.6cm} The assertion holds true similarly for both 
                the case of $ \alpha_j \leq \alpha'_j $, 
                $ \alpha_k > \alpha'_k $, \\
\hspace*{1.6cm} $ \alpha_l \leq \alpha'_l $, 
                and the case of $ \alpha_j > \alpha'_j $,
                $ \alpha_k \leq \alpha'_k $, $ \alpha_l \leq \alpha'_l $. \\  
\hspace*{0.6cm} ($\bold 1$-II-i-b) (No.86-88, 92) \ \ 
                $ z_i^I, \; z_j^J z_k, \; z_k^K z_l, \; z_j z_l^L \in f \quad
              (\exists I, \; \exists J, \; \exists K, \; \exists L \geq 2) $.\\
\hspace*{1.6cm} $ \alpha_j \leq \alpha'_j $ or $ \alpha_l \leq \alpha'_l $, 
                \ and $ \alpha_j+\alpha_l+(L-1)(\alpha_l-\alpha'_l) \leq 
                        \alpha'_j+\alpha'_l.$\\
\hspace*{1.6cm} The rest is similar to $\bold 1$-II-i-a. \\
\hspace*{0.6cm} ($\bold 1$-II-i-c) (No.84-89) \ \ 
                $ z_i^I, \; z_j^J z_k, \; z_k^K z_l, \; z_i z_l^L \in f \quad
              (\exists I, \; \exists J, \; \exists K, \; \exists L \geq 2) $.\\
\hspace*{1.6cm} $ \alpha_i+\alpha_l+(L-1)(\alpha_l-\alpha'_l) \leq 
                  \alpha'_i+\alpha'_l.$ \\
\hspace*{1.6cm} When $ \alpha_k \leq \alpha'_k $, the rest is similar to 
                $\bold 1$-I-iii-b. \\
\hspace*{1.6cm} If $ \alpha_k > \alpha'_k $ then 
                   $ \alpha_j \leq \alpha'_j $ and $ \alpha_l \leq \alpha'_l $,
                so \\
\hspace*{1.6cm} $ \alpha _i+\alpha _j+\alpha _k+\alpha _l 
                < \alpha _i+\alpha _j+\alpha _k+\alpha _l
                  +(K-1)(\alpha_k-\alpha'_k) \leq 
                  \alpha'_i+\alpha'_j+\alpha'_k+\alpha'_l . $ \\
\hspace*{0.4cm} ($\bold 1$-II-ii) \ \ 
                $ z_i^I, \; z_j^J z_k, \; z_j z_k^K \in f \quad
                  (\exists I, \; \exists J, \; \exists K \geq 2) $.\\
\hspace*{1.4cm} $ \alpha_j+\alpha_k+(K-1)(\alpha_k-\alpha'_k) \leq 
                  \alpha'_j+\alpha'_k.$\\
\hspace*{0.6cm} ($\bold 1$-II-ii-a) (No.89) \ \ 
                $ z_i^I, \; z_j^J z_k, \; z_j z_k^K, \; z_i z_l^L \in f \quad
               (\exists I, \; \exists J, \; \exists K, \; \exists L \geq 2)$.\\
\hspace*{1.6cm} $ \alpha_i+\alpha_l+(L-1)(\alpha_l-\alpha'_l) \leq 
                  \alpha'_i+\alpha'_l.$\\
\hspace*{1.6cm} For the case $ \alpha_l \leq \alpha'_l $, 
                the assertion holds true. \\
\hspace*{1.6cm} If $ \alpha_l > \alpha'_l $ then \\
\hspace*{1.6cm} $ \alpha_i+\alpha_j+\alpha_k+\alpha_l 
                < \alpha_i+\alpha_j+\alpha_k+\alpha_l+(L-1)(\alpha_l-\alpha'_l)
                $ \\
\hspace*{1.6cm} $ < \begin{cases}
                  & \alpha_i+\alpha_j+\alpha_k+\alpha_l
                    +(J-1)(\alpha_j-\alpha'_j)+(L-1)(\alpha_l-\alpha'_l) \quad
                    ( \text{if} \; \alpha_j>\alpha'_j ) \\
                  & \alpha_i+\alpha_j+\alpha_k+\alpha_l
                    +(K-1)(\alpha_k-\alpha'_k)+(L-1)(\alpha_l-\alpha'_l) \quad
                    ( \text{if} \; \alpha_k>\alpha'_k )
                  \end{cases} $ \\
\hspace*{1.6cm} $ \leq \alpha'_i+\alpha'_j+\alpha'_k+\alpha'_l $.\\
\hspace*{0.4cm} ($\bold 1$-II-iii) \ \ 
                $ z_i^I, \; z_j^J z_k, \; z_i z_k^K \in f \quad
                  (\exists I, \; \exists J, \; \exists K \geq 2)$.\\	
\hspace*{1.4cm} $ \alpha_i+\alpha_k+(K-1)(\alpha_k-\alpha'_k) \leq 
                  \alpha'_i+\alpha'_k.$\\
\hspace*{0.6cm} ($\bold 1$-II-iii-a) (No.89) \ \ 
                $ z_i^I, \; z_j^J z_k, \; z_i z_k^K, \; z_k z_l^L \in f \quad
               (\exists I, \; \exists J, \; \exists K, \; \exists L \geq 2)$.\\
\hspace*{1.6cm} $ \alpha_k \leq \alpha'_k $ or $ \alpha_l \leq \alpha'_l $, \ 
                and $ \alpha_k+\alpha_l+(L-1)(\alpha_l-\alpha'_l) \leq
                      \alpha'_k+\alpha'_l $.\\
\hspace*{1.6cm} Moreover 
                $ \alpha_j \leq \alpha'_j $ or $ \alpha_l \leq \alpha'_l $
                since there exists $ z_j^q z_l^s \in f $.\\
\hspace*{1.6cm} So $ \alpha_1+\dotsb+\alpha_4 < \alpha'_1+\dotsb+\alpha'_4 $ 
                similarly as $\bold 1$-II-i-a.\\
\hspace*{0.6cm} ($\bold 1$-II-iii-b) (No.85, 87, 89) \ \ 
                $ z_i^I, \; z_j^J z_k, \; z_i z_k^K, \; z_i z_l^L \in f \quad
               (\exists I, \; \exists J, \; \exists K, \; \exists L \geq 2)$.\\
\hspace*{1.6cm} Since there exists $ z_k^r z_l^s \in f $, we have 
                $ \alpha_k \leq \alpha'_k $ or $ \alpha_l \leq \alpha'_l $.\\
\hspace*{1.6cm} Moreover $ \alpha_i+\alpha_l+(L-1)(\alpha_l-\alpha'_l) \leq 
                           \alpha'_i+\alpha'_l $ by $ z_i z_l^L \in f $.\\
\hspace*{1.6cm} When $ \alpha_k \leq \alpha'_k $, the rest is similar to 
                $\bold 1$-I-iii-b. \\
\hspace*{1.6cm} If $ \alpha_k > \alpha'_k $ then 
                   $ \alpha_j \leq \alpha'_j $ and $ \alpha_l \leq \alpha'_l $,
                so \\
\hspace*{1.6cm} $ \alpha _i+\alpha _j+\alpha _k+\alpha _l 
                < \alpha _i+\alpha _j+\alpha _k+\alpha _l
                  +(K-1)(\alpha_k-\alpha'_k) \leq 
                  \alpha'_i+\alpha'_j+\alpha'_k+\alpha'_l . $ \\
Case $\bold 2.$ \ \ $ z_i^I z_j \in f \quad (\exists I \geq 2) $.\\
\hspace*{1cm} $ \alpha_i \leq \alpha'_i $ or $ \alpha_j \leq \alpha'_j $, \ 
              and $ \alpha_i+\alpha_j+(I-1)(\alpha_i-\alpha'_i) 
                    \leq \alpha'_i+\alpha'_j $.\\
\hspace*{0.2cm} ($\bold 2$-I) \ \ $ z_i^I z_j, \; z_j^J z_k \in f \quad
                                    (\exists I, \; \exists J \geq 2) $.\\
\hspace*{1.2cm} $ \alpha_j \leq \alpha'_j $ or $ \alpha_k \leq \alpha'_k $, \  
                and $ \alpha_j+\alpha_k+(J-1)(\alpha_j-\alpha'_j) 
                      \leq \alpha'_j+\alpha'_k $.\\                
\hspace*{0.4cm} ($\bold 2$-I-i) \ \ 
                $ z_i^I z_j, \; z_j^J z_k, \; z_k^K z_l \in f \quad
                  (\exists I, \; \exists J, \; \exists K \geq 2) $.\\
\hspace*{1.4cm} $ \alpha_k \leq \alpha'_k $ or $ \alpha_l \leq \alpha'_l $, \ 
                and $ \alpha_k+\alpha_l+(K-1)(\alpha_k-\alpha'_k) 
                      \leq \alpha'_k+\alpha'_l $.\\
\hspace*{0.6cm} ($\bold 2$-I-i-a) (No.94, 95) \ \ 
                $ z_i^I z_j, \; z_j^J z_k, \; z_k^K z_l, \; z_k z_l^L \in f 
                \quad 
               (\exists I, \; \exists J, \; \exists K, \; \exists L \geq 2)$.\\
\hspace*{1.6cm} $ \alpha_k+\alpha_l+(L-1)(\alpha_l-\alpha'_l) 
                  \leq \alpha'_k+\alpha'_l $.\\
\hspace*{1.6cm} Since there exists $ z_j^q z_l^s \in f $, we have 
                $ \alpha_j \leq \alpha'_j $ or $ \alpha_l \leq \alpha'_l $.\\
\hspace*{1.6cm} If $ \alpha_i \leq \alpha'_i $ then 
                the rest is similar to $\bold 1$-II-i-a. \\
\hspace*{1.6cm} If $ \alpha_i > \alpha'_i $ then $ \alpha_j \leq \alpha'_j $ 
                and so \\
\hspace*{1.6cm} $ \alpha_i+\alpha_j+\alpha_k+\alpha_l
                < \alpha_i+\alpha_j+\alpha_k+\alpha_l+(I-1)(\alpha_i-\alpha'_i)
                $ \\
\hspace*{1.6cm} $ < \begin{cases}
                  & \alpha_i+\alpha_j+\alpha_k+\alpha_l
                    +(I-1)(\alpha_i-\alpha'_i)+(K-1)(\alpha_k-\alpha'_k) \quad
                    ( \text{if} \; \alpha_k>\alpha'_k ) \\
                  & \alpha_i+\alpha_j+\alpha_k+\alpha_l
                    +(I-1)(\alpha_i-\alpha'_i)+(L-1)(\alpha_l-\alpha'_l) \quad
                    ( \text{if} \; \alpha_l>\alpha'_l )
                  \end{cases} $ \\
\hspace*{1.6cm} $ \leq \alpha'_i+\alpha'_j+\alpha'_k+\alpha'_l.$ \\
\hspace*{0.6cm} ($\bold 2$-I-i-b) (No.94, 95) \ \ 
                $ z_i^I z_j, \; z_j^J z_k, \; z_k^K z_l, \; z_j z_l^L \in f 
                \quad 
               (\exists I, \; \exists J, \; \exists K, \; \exists L \geq 2)$.\\
\hspace*{1.6cm} $ \alpha_j \leq \alpha'_j $ or $ \alpha_l \leq \alpha'_l $, \ 
                and $ \alpha_j+\alpha_l+(L-1)(\alpha_l-\alpha'_l)
                      \leq \alpha'_j+\alpha'_l $.\\
\hspace*{1.6cm} Since there exists $ z_i^p z_l^s \in f $, we have 
                $ \alpha_i \leq \alpha'_i $ or $ \alpha_l \leq \alpha'_l $.\\
\hspace*{1.6cm} If $ \alpha_i \leq \alpha'_i $ then 
                the rest is similar to $\bold 1$-II-i-a. \\
\hspace*{1.6cm} If $ \alpha_i > \alpha'_i $ then 
                $ \alpha_j \leq \alpha'_j $ and $ \alpha_l \leq \alpha'_l $, \ 
                and so \\
\hspace*{1.6cm} $ \alpha_i+\alpha_j+\alpha_k+\alpha_l 
                < \alpha_i+\alpha_j+\alpha_k+\alpha_l+(I-1)(\alpha_i-\alpha'_i)
                  +(K-1)(\alpha_k-\alpha'_k) $ \\
\hspace*{4.8cm} $ \leq \alpha'_i+\alpha'_j+\alpha'_k+\alpha'_l $ \qquad
                for the case of $ \alpha_k > \alpha'_k $.\\
\hspace*{0.6cm} ($\bold 2$-I-i-c) (No.94, 95) \ \ 
                $ z_i^I z_j, \; z_j^J z_k, \; z_k^K z_l, \; z_i z_l^L \in f 
                \quad 
               (\exists I, \; \exists J, \; \exists K, \; \exists L \geq 2)$.\\
\hspace*{1.6cm} $ \alpha_i \leq \alpha'_i $ or $ \alpha_l \leq \alpha'_l $, \ 
                and $ \alpha_i+\alpha_l+(L-1)(\alpha_l-\alpha'_l)
                      \leq \alpha'_i+\alpha'_l $.\\
\hspace*{1.6cm} If $ \alpha_i \leq \alpha'_i $ and $ \alpha_j \leq \alpha'_j $ 
                then the assertion holds true. \\
\hspace*{1.6cm} If $ \alpha_i \leq \alpha'_i $ and $ \alpha_j > \alpha'_j $ 
                then $ \alpha_k \leq \alpha'_k $, so the rest is similar to 
                $\bold 1$-I-iii-b. \\
\hspace*{1.6cm} If $ \alpha_i > \alpha'_i $ and $ \alpha_j \leq \alpha'_j $ 
                then $ \alpha_l \leq \alpha'_l $, so the rest is similar to 
                $\bold 1$-I-iii-b. \\
\hspace*{0.2cm} ($\bold 2$-II) \ \ $ z_i^I z_j, \; z_i z_j^J \in f \quad 
                                     (\exists I, \; \exists J \geq 2) $.\\
\hspace*{1.2cm} $ \alpha_i+\alpha_j+(J-1)(\alpha_j-\alpha'_j) 
                  \leq \alpha'_i+\alpha'_j $.\\
\hspace*{0.4cm} ($\bold 2$-II-i) \ \ 
                $ z_i^I z_j, \; z_i z_j^J, \; z_k^K z_l \in f \quad 
                  (\exists I, \; \exists J, \; \exists K \geq 2) $.\\
\hspace*{1.4cm} $ \alpha_k \leq \alpha'_k $ or $ \alpha_l \leq \alpha'_l $, \ 
                and $ \alpha_k+\alpha_l+(K-1)(\alpha_k-\alpha'_k)
                      \leq \alpha'_k+\alpha'_l $.\\
\hspace*{0.6cm} ($\bold 2$-II-i-a) \ \ 
                $ z_i^I z_j, \; z_i z_j^J, \; z_k^K z_l, \; z_k z_l^L \in f 
                \quad 
               (\exists I, \; \exists J, \; \exists K, \; \exists L \geq 2)$.\\
\hspace*{1.6cm} $ \alpha_k+\alpha_l+(L-1)(\alpha_l-\alpha'_l)
                  \leq \alpha'_k+\alpha'_l $.\\
\hspace*{1.6cm} $ \alpha_i+\alpha_j+\alpha_k+\alpha_l $ \\
\hspace*{1.6cm} $ < \begin{cases}
                  & \alpha_i+\dotsb+\alpha_l+(I-1)(\alpha_i-\alpha'_i) \qquad
                  ( \text{if} \; \alpha_i>\alpha'_i ) \\
                  & \dotsb\dotsb\dotsb\dotsb \hspace{4.5cm} \dotsb\dotsb \\
                  & \alpha_i+\dotsb+\alpha_l+(L-1)(\alpha_l-\alpha'_l) \qquad
                  ( \text{if} \; \alpha_l>\alpha'_l ) 
                  \end{cases} $ \\
\hspace*{1.6cm} $ < \begin{cases}
                  & \alpha_i+\dotsb+\alpha_l
                    +(I-1)(\alpha_i-\alpha'_i)+(K-1)(\alpha_k-\alpha'_k) \quad
                  ( \alpha_i>\alpha'_i, \; \alpha_k>\alpha'_k ) \\
                  & \alpha_i+\dotsb+\alpha_l
                    +(I-1)(\alpha_i-\alpha'_i)+(L-1)(\alpha_l-\alpha'_l) \quad
                  ( \alpha_i>\alpha'_i, \; \alpha_l>\alpha'_l ) \\
                  & \alpha_i+\dotsb+\alpha_l
                    +(J-1)(\alpha_j-\alpha'_j)+(K-1)(\alpha_k-\alpha'_k) \quad
                  ( \alpha_j>\alpha'_j, \; \alpha_k>\alpha'_k ) \\
                  & \alpha_i+\dotsb+\alpha_l
                    +(J-1)(\alpha_j-\alpha'_j)+(L-1)(\alpha_l-\alpha'_l) \quad
                  ( \alpha_j>\alpha'_j, \; \alpha_l>\alpha'_l )
                  \end{cases} $ \\
\hspace*{1.6cm} $ \leq \alpha'_i+\alpha'_j+\alpha'_k+\alpha'_l.$ \\
\hspace*{0.4cm} ($\bold 2$-II-ii) \ \ $ z_i^I z_j, \; z_i z_j^J, \; z_j z_k^K 
                                        \in f \quad
                            (\exists I, \; \exists J, \; \exists K \geq 2) $.\\
\hspace*{1.4cm} $ \alpha_j \leq \alpha'_j $ or $ \alpha_k \leq \alpha'_k $, \ 
                and $ \alpha_j+\alpha_k+(K-1)(\alpha_k-\alpha'_k)
                      \leq \alpha'_j+\alpha'_k $.\\
\hspace*{0.6cm} ($\bold 2$-II-ii-a) \ \ $ z_i^I z_j, \; z_i z_j^J, \; 
                                          z_j z_k^K, \; z_j z_l^L \in f \quad 
              (\exists I, \; \exists J, \; \exists K, \; \exists L \geq 2) $.\\
\hspace*{1.6cm} $ \alpha_j \leq \alpha'_j $ or $ \alpha_l \leq \alpha'_l $, \ 
                and $ \alpha_j+\alpha_l+(L-1)(\alpha_l-\alpha'_l)
                      \leq \alpha'_j+\alpha'_l $.\\
\hspace*{0.6cm} ($\bold 2$-II-ii-b) (No.94, 95) \ \ 
                $ z_i^I z_j, \; z_i z_j^J, \; z_j z_k^K, \; z_i z_l^L \in f 
                \quad 
              (\exists I, \; \exists J, \; \exists K, \; \exists L \geq 2) $.\\
\hspace*{1.6cm} $ \alpha_i \leq \alpha'_i $ or $ \alpha_l \leq \alpha'_l $, \ 
                and $ \alpha_i+\alpha_l+(L-1)(\alpha_l-\alpha'_l) 
                      \leq \alpha'_i+\alpha'_l $.\\
\hspace*{1.6cm} Since there exists $ z_i^p z_k^r \in f $, we have 
                $ \alpha_i \leq \alpha'_i $ or $ \alpha_k \leq \alpha'_k $.\\
\hspace*{1.6cm} Since there exists $ z_j^q z_l^s \in f $, we have
                $ \alpha_j \leq \alpha'_j $ or $ \alpha_l \leq \alpha'_l $.\\
\hspace*{1.6cm} If $ \alpha_l \leq \alpha'_l $ then the rest is similar to 
                $\bold 1$-II-i-a. \\
\hspace*{1.6cm} If $ \alpha_l > \alpha'_l $ then $ \alpha_i \leq \alpha'_i $ 
                and $ \alpha_j \leq \alpha'_j $, so the rest is similar to 
                $\bold 2$-I-i-b. \\
\end{pf} 
\begin{lem}
Let $ f=\alpha(x,y,z)x^2+\beta(y,z)y^3+\phi(z)y^2+\varphi(z)y+\psi(z) \in 
{\bold C}[x,y,z] $ define an isolated singularity at the origin $ 0 \in 
{\bold C}^3 $, which satisfies one of the following\rom:

  $(1)$ \ $ \beta(y,z)=\beta_0+\text{higher terms}, \; 
            0 \ne \beta_0 \in {\bold C}, \; \text{ord}\phi=1 $, and \\
\hspace*{1.1cm} $ \beta_0 y^3+\phi_0(z)y^2+\varphi_0(z)y+\psi_0(z) $
                has no triple factor, 

  $(1')$ \ $ \beta \equiv 0 $ or $ \text{ord}\beta \geq 1 $, \ and 
           $ \text{ord}\phi=1 $,

  $(2)$ \ $ \beta(y,z)=\beta_0+\text{higher terms}, \; 
            0 \ne \beta_0 \in {\bold C}, \; \;  
            \text{ord}\varphi \leq 3 $ or $ \text{ord}\psi \leq 5 $, \ and \\
\hspace*{1.1cm} $ \beta_0 y^3+\phi_0(z)y^2+\varphi_0(z)y+\psi_0(z) $ 
                has no triple factor, \\
where $ \text{ord}\alpha=0 $, and 
$ \phi_0, \; \varphi_0, \; \psi_0 $ are the initial parts of 
$ \phi, \; \varphi, \; \psi $, respectively. 
Then $(f,0)$ is rational.
\end{lem}
\begin{pf} 
\hspace*{0.2cm} We may assume that $ \alpha(x,y,z)=1, \; \beta_0=1 $. \\
(1) \ If $ \text{ord}\varphi \leq 1 $ or $ \text{ord}\psi \leq 2 $ then the 
assertion holds true. So we may assume $ \text{ord}\varphi \geq 2 $ and 
$ \text{ord}\psi \geq 3 $. Let $f_0$ be the initial part of $f$ with respect 
to the weight $ \alpha = (1/2, \; p, \; q) $, where $\alpha$ satisfies the 
conitions 
$$ \frac13 \leq p = \frac{1-q}{2} < \frac12 < p+q \left( = \frac{1+q}{2} 
   \right), $$
$$ \deg_{\alpha}(\varphi(z)y), \; \deg_{\alpha}\psi(z) \geq 1, \quad 
   \text{and} \quad 
   \deg_{\alpha}(\varphi(z)y) \; \text{or} \; \deg_{\alpha}\psi(z)=1. $$
(There exists such weight $\alpha$ because 
$$ \frac{1-q}{2} \geq \frac13, \; \frac{1-q}{2}+Nq \geq 1 \quad  
\text{for  } N \geq 2, \; \frac13 \geq q \geq \frac{1}{2N-1}.) $$
If $ f_0-x^2 \in {\bold C}[y,z] $ has no double factor, then $(f_0, 0)$ is an 
isolated singularity, so $(f_0, 0)$ is rational. So $ \{f_0=0\} \subset 
{\bold C}^3 $ has only rational singularities around the origin $ 0 \in
{\bold C}^3 $. Therefore $(f,0)$ is also rational from 2.4 and 2.5, because 
$ \displaystyle \frac12+p+q>1 $.\\
If $ f_0-x^2 $ has a double factor, namely, \\
\hspace*{1cm} $ f=x^2+(y+\gamma_1 z)^2 (y+\gamma_2 z)+\text{higher terms}, \; 
\gamma_1 \ne \gamma_2 $, \\
then taking the coordinate changes $ Y:=y+\gamma_1 z $ and 
                                   $ z':=(\gamma_2-\gamma_1)z $, thus \\
\hspace*{1cm} $ f=x^2+\beta'(Y,z')Y^3+\phi'(z')Y^2+\varphi'(z')Y+\psi'(z')$,\\ 
for some $ \beta' \in {\bold C}[Y,z'], \; 
           \phi', \; \varphi', \; \psi' \in {\bold C}[z'] $ with  
$ \beta'=1+\text{higher terms}, \; 
  \phi'=z'+\text{higher terms}, \; 
  2=\text{ord}\varphi<\text{ord}\varphi', \; 
  3=\text{ord}\psi<\text{ord}\psi' $. 
Let $ \varphi'_0, \; \psi'_0 $ be the initial parts of $ \varphi', \; 
\psi' $, respectively. We replace the weight $ \alpha = (1/2, \; p, \; q) $
with $ \alpha' = (1/2, \; p', \; q') $, which satisfies: 
$$ \frac13 = p < p' = \frac{1-q'}{2} < \frac12 
   < p'+q' \left( = \frac{1+q'}{2} \right), $$
$$ \deg_{\alpha'}(\varphi'(z')Y), \; \deg_{\alpha'}\psi'(z') \geq 1, 
   \quad \text{and} \quad
   \deg_{\alpha'}(\varphi'(z')Y) \; \text{or} \; 
   \deg_{\alpha'}\psi'(z')=1. $$
(There exists such weight $\alpha'$ because 
$$ \frac{1-q'}{2} > \frac13, \; \frac{1-q'}{2}+Nq' \geq 1 \quad  
\text{for  } N > 2, \; \frac13 > q' \geq \frac{1}{2N-1}.) $$
If \ $ Y^2 z'+Y\varphi'_0(z')+\psi'_0(z')=(Y+g(z'))^2 z' $ \; 
for some $ g \in {\bold C}[z'] \; \; (\text{ord}(g) > 1) $,\\
then taking the coordinate change $ Y':=Y+g(z') $, thus \\
\hspace*{1cm} $ f=x^2+\beta''(Y',z'){Y'}^3+\phi''(z'){Y'}^2
                  +\varphi''(z')Y'+\psi''(z') $, \\
for some $ \beta'' \in {\bold C}[Y',z'], \; 
           \phi'', \; \varphi'', \; \psi'' \in {\bold C}[z'] $ with
$ \beta''=1+\text{higher terms}, \; 
  \phi''=z'+\text{higher terms}, \; 
  \text{ord} \varphi' < \text{ord} \varphi'', \; 
  \text{ord} \psi' < \text{ord} \psi'' $. 
Let $ \varphi''_0, \; \psi''_0 $ be the initial parts of 
$ \varphi'', \; \psi'' $, respectively. If 
$$ {Y'}^2 z'+Y'\varphi''_0(z')+\psi''_0(z')=(Y'+g'(z'))^2 z' $$ 
for some $ g' \in {\bold C}[z'] \; (1 < \text{ord}(g) < \text{ord}(g')) $, 
then taking the coordinate change \\ $ Y'':=Y'+g'(z') $, $\dotsb\dotsb$. 
If this procedure continues infinitely, then \\
\hspace*{2.5cm} $ \displaystyle 
                \frac13=p<p'<\dotsb<p^{(n)}=\frac{1-q^{(n)}}{2}<\dotsb
                <\frac12<p^{(n)}+q^{(n)} $,\\
\hspace*{2.5cm} $ \displaystyle \frac13=q>q'>\dotsb>q^{(n)}>\dotsb>0 $,
$$ \dim_{\bold C}{\bold C}\{z'\} \left/
   \left(\varphi^{(n)}, \frac{\partial \psi^{(n)}}
  {\partial z'} \right) \right. \leq \mu(f,0) 
   \overset{n \to +\infty}{\longrightarrow} +\infty, $$ 
a contradiction. Therefore $(f,0)$ is rational by 2.4 and 2.5. \\
(1$'$) \ We may assume $ \phi=z+\text{higher terms}, \; 
\text{ord}\varphi \geq 2 $ and $ \text{ord}\psi \geq 3 $. 
Let $f_0$ be the initial part of $f$ with respect to the weight 
$ \alpha = (1/2, \; p, \; q) $, where $\alpha$ satisfies: 
$$ \frac13 \leq p = \frac{1-q}{2} < \frac12 < p+q \left( = \frac{1+q}{2} 
   \right), $$
$$ \deg_{\alpha}(\varphi(z)y), \; \deg_{\alpha}\psi(z) \geq 1, \quad 
   \text{and} \quad 
   \deg_{\alpha}(\varphi(z)y) \; \text{or} \; \deg_{\alpha}\psi(z)=1. $$
If $ f_0-x^2 \in {\bold C}[y,z] $ has no double factor, then $(f_0, 0)$ is 
rational similarly as in (1). \\
If $ f_0-x^2 $ has a double factor, namely, \\
\hspace*{1cm} $ y^2 z+y\varphi_0(z)+\psi_0(z)=(y+g(z))^2 z $ \; 
              for some $ g \in {\bold C}[z] \; \; (\text{ord}(g) \geq 1) $,\\
(where $ \varphi_0, \; \psi_0 $ are the initial parts of $ \varphi, \; 
\psi $, respectively,) 
then taking the coordinate change $ Y:=y+g(z) $, thus \\
\hspace*{1cm} $ f=
  \begin{cases}
     x^2+\phi(z)Y^2+\varphi'(z)Y+\psi'(z), \qquad & (\beta \equiv 0) \\
     x^2+\beta'(Y,z)Y^3+\phi'(z)Y^2+\varphi'(z)Y+\psi'(z), \qquad  
     & (\text{ord}\beta \geq 1) 
  \end{cases} $ \\
for some $ \beta' \in {\bold C}[Y,z], \; 
           \phi', \; \varphi', \; \psi' \in {\bold C}[z] $ with  
$ 1 \leq \text{ord}\beta', \; \phi'=z+\text{higher terms}, \; 
  2 \leq \text{ord} \varphi < \text{ord} \varphi', \; 
  3 \leq \text{ord} \psi < \text{ord} \psi' $. 
We replace the weight $ \alpha = (1/2, \; p, \; q) $
with $ \alpha' = (1/2, \; p', \; q') $, which satisfies the conitions
$$ \frac13 \leq p < p' = \frac{1-q'}{2} < \frac12 
   < p'+q' \left( = \frac{1+q'}{2} \right), $$
$$ \deg_{\alpha'}(\varphi'(z)Y), \; \deg_{\alpha'}\psi'(z) \geq 1, \quad
   \text{and} \quad
   \deg_{\alpha'}(\varphi'(z)Y) \; \text{or} \; 
   \deg_{\alpha'}\psi'(z)=1. $$
If this procedure continues infinitely, then 
$ \mu(f,0) \overset{n \to +\infty}{\longrightarrow} +\infty, $ 
a contradiction. \\ 
(2) \ Let $f_0$ be the initial part of $f$ with respect to the weight 
$ \alpha = (1/2, \; 1/3, \; q) $.\\
If $ f_0-x^2 \in {\bold C}[y,z] $ has no double factor, then $(f_0, 0)$ is an 
isolated singularity, and so $(f_0, 0)$ is rational for the case of 
$ \text{ord}\varphi \leq 3 $ or $ \text{ord}\psi \leq 5 $. Therefore $(f,0)$ is
also rational from 2.4 and 2.5, because 
$ \displaystyle \frac12+\frac13+q>1 $ for $ \displaystyle q>\frac16 $.\\
So consider the case of 
$ f=x^2+(y+\gamma_1 z)^2 (y+\gamma_2 z)+\text{higher terms}, \; 
\gamma_1 \ne \gamma_2 $. \\
Then the situation is similar to (1). 
\end{pf}
\begin{exmp} 
\hspace*{0.3cm} (See Theorem 2.7.) \\
\hspace*{1cm} $ f = x^2 + y^3 + z^9 + w^{18} $, \ \ \ (No.12) \\
\hspace*{1cm} $ f_{\lambda} = x^2 + y^3 + 
                (- \frac{27}{4}c_{3 0}^2 (z+\gamma_1 w)^4 (z-2\gamma_1 w)^2
                 + \varphi(z,w))y $ \\
\hspace*{2cm} $ + (z+\gamma_1 w)^6 (z-2\gamma_1 w)^3 + \psi(z,w) $,\\
where $ - \frac{27}{4}c_{3 0}^3 = 1, \; \text{ord}\varphi \geq 7, \; 
        \text{ord}\psi \geq 10, \; \psi(z,w) \ni w^{18} $ \ (Case 10-b). \\
Taking the coordinate change 
$ Y:= y + \frac32 c_{3 0} (z+\gamma_1 w)^2 (z-2\gamma_1 w) $, we have \\
\hspace*{1cm} $ f_{\lambda} = x^2 + Y^3 
                - \frac92 c_{3 0} (z+\gamma_1 w)^2 (z-2\gamma_1 w)Y^2 
                + \varphi(z,w)Y $ \\
\hspace*{2cm} $ - \frac32 c_{3 0} (z+\gamma_1 w)^2 (z-2\gamma_1 w) \varphi(z,w)
                + \psi(z,w) $,\\
and then $ z':= z+\gamma_1 w $,\\
\hspace*{1cm} $ f_{\lambda} = x^2 + Y^3
                - \frac92 c_{3 0} ({z'}^3 - 3\gamma_1 {z'}^2 w)Y^2 
                + \varphi'(z',w)Y + \psi'(z',w) $ \\
for some $ \varphi', \; \psi' \in {\bold C}[z',w], \; 
\text{ord}\varphi' \geq 7, \; \text{ord}\psi' \geq 10 $.\\
Let $ \alpha':=(1/2, 1/3, 1/8, 1/12), \; \alpha'':=(1/2, 17/50, 3/25, 2/25), 
\; c'_{2 1}:=-3\gamma_1 c_{3 0} $, and \\
\hspace*{1cm} $ \varphi'=-3(3c'_{2 1}{z'}^2 w+c'_{0 4}w^4)c'_{0 4}w^4+\Phi', 
                \qquad \quad \; \Phi'_0=c'_{1 7}z'w^7 $,\\
\hspace*{1cm} $ \psi'=(-\frac92c'_{2 1}{z'}^2 w-2c'_{0 4}w^4)(c'_{0 4}w^4)^2
                       +\Psi', \qquad 
                \Psi'_0=c'_{5 5}{z'}^5 w^5+c'_{1 \; 11}z' w^{11} $,\\
where $ \Phi'_0, \; \Psi'_0 $ are the initial parts of $ \Phi', \; \Psi' $ 
with respect to the weight $ \alpha' $. Then \\
\hspace*{1cm} $ f_{\lambda} = x^2 + Y^3 
                - \frac92 (c_{3 0}{z'}^3+c'_{2 1}{z'}^2 w)Y^2 + 
                (-3(3c'_{2 1}{z'}^2 w+c'_{0 4}w^4)c'_{0 4}w^4+\Phi')(z',w)Y $\\
\hspace*{2cm} $ +((-\frac92c'_{2 1}{z'}^2 w-2c'_{0 4}w^4)(c'_{0 4}w^4)^2+\Psi')
                 (z',w) $ \\
\hspace*{1.5cm} $ = x^2+(Y+c'_{0 4}w^4)^2
                        (Y-\frac92c'_{2 1}{z'}^2 w-2c'_{0 4}w^4) $ \\
\hspace*{2cm} $     -\frac92c_{3 0}{z'}^3 Y^2+\Phi'(z',w)Y+\Psi'(z',w) $ 
                \qquad \qquad (Case II-A-II$'$-ii-i). \\
Taking the coordinate change $ Y':=Y+c'_{0 4}w^4 $, we have \\
\hspace*{1cm} $ f_{\lambda} = x^2 + {Y'}^3 - (\frac92c_{3 0}{z'}^3
                              +\frac92c'_{2 1}{z'}^2 w+3c'_{0 4}w^4){Y'}^2 $ \\
\hspace*{2cm} $               + (9c_{3 0}c'_{0 4}{z'}^3 w^4+\Phi'(z',w))Y' 
                              - \frac92c_{3 0}{c'}_{0 4}^2 {z'}^3 w^8 
                              - c'_{0 4}w^4\Phi'(z',w) + \Psi'(z',w) $,\\
\hspace*{1cm} $ f_0 = x^2 - (\frac92c'_{2 1}{z'}^2 w+3c'_{0 4}w^4){Y'}^2 
                      - \frac92c_{3 0}{c'}_{0 4}^2 {z'}^3 w^8 
                      - c'_{0 4}w^4\Phi'_0(z',w) + \Psi'_0(z',w) $.\\
Then $(f_0, 0)$ is an isolated singularity or Case II-A-I, 
and $ 1/2+17/50+3/25+2/25>1 $, so $ \{f_0=0\}-\{0\} \subset {\bold C}^4 $ has 
at most rational singularities around the origin $ 0 \in {\bold C}^4 $. 
\end{exmp}
\begin{rem}

  We will take the following arguments to show Theorem 2.7. 

  Taking suitable local coordinate changes finitely, if necessary, we have a 
three-dimensional face $\Delta$ of $\Gamma(f_{\lambda})$ such that $(1,1,1,1)$
lies strictly above the hyperplane including $\Delta$, that 
$ \operatorname{Sing}(f_0) := \operatorname{Sing}(\{f_0=0\}) = \bigcup_j C_j $
with $ \dim_{\bold C} C_j \leq 1 $, and that $ (f_0, P):=(\{f_0=0\}, P) $ is 
rational for any $ P \in C_j \subset \operatorname{Sing}(f_0) $ with 
$ P \ne 0 $ and $ 0 \in C_j $, 
where $ f_0 := f_{\lambda \; \Delta} = \sum_{\nu \in \Delta} a_{\nu} z^{\nu} 
\subset f_{\lambda} $. Therefore $ \{f_0=0\}-\{0\} $ has at most rational 
singularities around $0$, and so $(f_{\lambda},0)$ is rational. 

  More precisely, any three-dimensional face 
$ \Delta = \{ (X,Y,Z,W) \in {{\bold R}_{\geq 0}}^4 \; ; \; 
\frac12 X + \beta Y + \gamma Z + \delta W = 1, \; 
0 \leq X \leq 2, \; j \leq Y \leq j', \; k \leq Z \leq k', \; l \leq W \leq l' 
\} $ of $\Gamma(f_{\lambda})$ satisfies $ 1 < \frac12 + \beta + \gamma + 
\delta $ under the original local coordinates from Lemma 3.1. 
Choose a three-dimensional face $\Delta$ of $\Gamma(f_{\lambda})$ suitably and 
let $ f_0 := f_{\lambda \; \Delta} = \sum_{\nu \in \Delta} a_{\nu} z^{\nu} 
\subset f_{\lambda} := \sum_{\nu} a_{\nu} z^{\nu} $.\\
\\
Case (0). \ $ (f_0,0) $ is an isolated singularity. 

  Then $ (f_0,0) $ is rational and so $ \{ f_0=0 \} \subset {\bold C}^4 $ has 
only rational singularities around the origin $ 0 \in {\bold C}^4 $. \\
\\
Case (I). \ $ \operatorname{Sing}(f_0) = \bigcup_j C_j, \qquad 
              \dim_{\bold C} C_j \leq 1 $ for all j. 

  Since $f_0$ is quasi-homogeneous, we have $ \operatorname{Sing}(C_j)=\{0\} $.
If an arbitrary point $ P \in C_j-\{0\} $ is rational on $\{ f_0=0 \}$ for each
irreducible curve $C_j$ with $ 0 \in C_j $, then $ \{ f_0=0 \}-\{ 0 \} $ has at
most rational singularities around 0. 
Rationality of $P$ is shown by using the following fact: \\
\hspace*{1cm} if there exists an hyperplane cut $ (H,P) \subset 
              ( \{ f_0=0 \}, \; P ) $ which is rational \\
\hspace*{1cm} and Gorenstein, then $ ( \{ f_0=0 \}, \; P ) $ is also rational 
              and Gorenstein.

  When we can not tell whether $P$ is rational or not, take suitable local 
coordinate change and repeat the same procedure as above. The condition
$ 1 < \frac12 + \beta^{(n)} + \gamma^{(n)} + \delta^{(n)} $ is still satisfied 
for a certain $\Delta^{(n)}$ after each coordinate change. This procedure must 
finish in finite times from the assumption $ \mu(f_{\lambda},0)<\mu(f,0) $.\\
\\
Case (II). \ $ \operatorname{Sing}(f_0) = \bigcup_j C_j, \qquad 
               \dim_{\bold C} C_j = 2 $ for some j. 

  After suitable local coordinate change, choose another three-dimensional 
face $\Delta$ of $\Gamma(f_{\lambda})$ properly, and 
let $ f_0 := f_{\lambda \; \Delta} = \sum_{\nu \in \Delta} a_{\nu} z^{\nu} $.

  (II-I$^{(\prime)}$). \ $ \operatorname{Sing}(f_0) = \bigcup_j C_j, \qquad    
                           \dim_{\bold C} C_j \leq 1 $ for all j. \\
\hspace*{0.5cm} Then the proof is completed similarly as in (I). 

  (II-II$^{(\prime)}$). \ $ \operatorname{Sing}(f_0) = \bigcup_j C_j, \qquad 
                            \dim_{\bold C} C_j = 2 $ for some j. \\
\hspace*{0.5cm} After suitable local coordinates change, the condition 
$ 1 < \frac12 + \beta' + \gamma' + \delta' $ is still satisfied 
for a certain face $ \Delta'= \{ (X,Y,Z,W) \; ; \;
\frac12 X + \beta' Y + \gamma' Z + \delta' W = 1 \} $ of $\Gamma(f_{\lambda})$.
Choose such three-dimensional face $ \Delta'$ properly; 
if Case (I$^{(\prime)}$) then the assertion is \\ concluded. 
If Case (II$^{(\prime)}$) again, take suitable coordinate change once more. 
This \\ procedure must finish in finite times from the assumption 
$ \mu(f_{\lambda},0)<\mu(f,0) $.
\end{rem}
{\it Proof of Theorem 2.7}.

  If there exists $ y^j z^k w^l \in f_{\lambda} $ such that $ j+k+l \leq 2 $ 
then $(f_{\lambda},0)$ is at most rational. So we may assume $ j+k+l \geq 3 $ 
for all $ y^j z^k w^l \in f_{\lambda} $.

  Let $ \Lambda :=\{(k,l) \; ; z^k w^l \in f_{\lambda}\} \bigcup\{\frac32(k,l) 
\; ; y z^k w^l \in f_{\lambda}\} $ and let $\Gamma$ be the union of the 
compact faces of the convex hull of $ \bigcup_{\nu \in \Lambda} (\nu + 
{{\bold R}_{\geq 0}}^2) $ in $ {\bold R}^2 $. For any one-dimensional face 
$ \Delta = \{ (Z,W) \; ; \; \gamma Z + \delta W = 1 \; , k_1 \leq Z \leq k_0, 
\; l_0 \leq W \leq l_1 \} $ of $\Gamma$, we get $ \gamma+\delta>1/6 $ from 
Lemma 3.1. Choose such $\Delta$ satisfies $ k_1<6 $ and $ l_0<6 $. (See 
{\sc Figure} 1.) Then, \\
\hspace*{1cm} $ f_{\lambda} = x^2+y^3+ \sum_{\gamma k+\delta l \geq 2/3} b_{kl}
                y z^k w^l + \sum_{\gamma k+\delta l \geq 1} c_{kl} z^k w^l $,\\
\hspace*{1cm} $ f_0 = x^2+ y^3+ \sum_{\gamma k+\delta l = 2/3} b_{kl} y z^k w^l
                + \sum_{\gamma k+\delta l = 1} c_{kl} z^k w^l $ \\
with respect to the weight $ \alpha:=(1/2, 1/3, \gamma, \delta) $.\\
\begin{figure}[h]
\begin{center}
\setlength{\unitlength}{1mm}
\begin{picture}(85,85)(-15,-15)
  \put(0,0){\vector(1,0){60}}
  \put(0,0){\vector(0,1){60}}
  \put(-5,-5){0}
  \put(63,-2){$Z$}
  \put(-2,63){$W$}
  \put(55,5){\line(-2,1){20}}
  \put(35,15){\thicklines\line(-1,1){20}}
  \put(15,35){\line(-1,4){5}}
  \put(30,30){\circle*{1}}
  \put(30,30){\makebox(15,6){(6,6)}}
  \put(35,15){\makebox(17,6){($k_0,l_0$)}}
  \put(15,35){\makebox(17,6){($k_1,l_1$)}}
  \put(15,50){$\Gamma$}
  \multiput(50,0)(-1,1){15}{\circle*{0.2}}
  \multiput(15,35)(-1,1){15}{\circle*{0.2}}
  \put(50,-10){$ \displaystyle \frac{1}{\gamma} $}
  \put(-5,50){$ \displaystyle \frac{1}{\delta} $}
  \put(20,20){$\Delta$}
\end{picture}
\caption{ }
\end{center}
\end{figure} \\
Case I. \ $ h:=f_0-x^2 \in {\bold C} [y,z,w] $ \ has no double factor. \\
Since $ h $ and $ \displaystyle \frac{\partial h}{\partial y} $ have no common 
factor, $ \dim_{\bold C} \operatorname{Sing}(f_0) \leq 1 $. Let $ C_j \subset 
\operatorname{Sing}(f_0) $ be an irreducible curve with $ 0 \in C_j $. 
If $ C_j \ni P = (0,a,b,c) \ne (0,0,0,0) $ then $ b \ne 0 $ or $ c \ne 0 $.\\
Let $ P \ne 0 $ be an arbitrary point on $C_j$. \\
\\
(I-A). \ $ C_j \ni P = (0,a(t),b(t),t) \; ; \; t \ne 0 $.\\
Since $f_0$ is quasi-homogeneous, we have $ a(t)=a't^{1/3\delta}, \; 
b(t)=b't^{\gamma/\delta} $ for some $ a', \; b' \in {\bold C} $.\\
Let $ \eta:=y-a, \; \zeta:=z-b $, \ and $ f_0(w=t):=f_0(x, \eta+a, \zeta+b, t) 
$.\\
Then $ f_0(w=t)-x^2 = {\eta}^3+3a{\eta}^2+\dotsb \in {\bold C} [\eta,\zeta] $ 
has no double factor. \\
(In fact, \ if $ f_0(w=t)-x^2=(\eta+\varphi(\zeta))^2 (\eta+\psi(\zeta)) $ then
{\allowdisplaybreaks
\begin{align*}
 f_0(w=t) - x^2 & = {\eta}^3 + 3a{\eta}^2 - 3 \varphi(\varphi-2a)\eta 
                    - {\varphi}^2 (2\varphi-3a) \\
                & = y^3 - 3(\varphi(z-b(t))-a(t))^2 y 
                    - 2(\varphi(z-b(t))-a(t))^3 ,\\
 f_0 - x^2 & = (y+\varphi(z-b(w))-a(w))^2 (y-2\varphi(z-b(w))+2a(w)) \\
           & = y^3 - 3(\varphi(z-b(w))-a(w))^2 y - 2(\varphi(z-b(w))-a(w))^3 .
\end{align*}}
It follows that $ \varphi(z-b(w))-a(w) \in {\bold C} [z,w] $ from 
$ \varphi(z-b(w))-a(w) \in {\bold C} [z,w^{\gamma/\delta}] $,\\ 
$(\varphi(z-b(w))-a(w))^2, \; (\varphi(z-b(w))-a(w))^3 \in {\bold C} [z,w]$.)\\
Therefore $ (f_0(w=t), \; (0,0,0)) $ is an isolated singularity under the local
coordinate system $ (x,\eta,\zeta) $.\\
If $ a \ne 0 $ then $(f_0(x,y,z,t), \; (0,a,b))$ is rational, so 
$(f_0, \; (0,a,b,t))$ is rational.\\
If $ a=0 $ then $ \eta=y $,  $ f_0(w=t)=x^2+y^3+\sum b_{k l} y (\zeta+b)^k t^l
+\sum c_{k l} (\zeta+b)^k t^l $.\\
If $ b=0 $ then $ \zeta=z $, $ f_0(w=t)=x^2+y^3+\sum b_{k l} y z^k t^l 
+\sum c_{k l} z^k t^l $. By Lemma 3.2, $ (f_0(x,y,z,t), \; (0,0,0)) $ is 
rational, so $ (f_0, \; (0,0,0,t)) $ is rational. \\
So we consider the case $ b \ne 0 $.\\
If there exists $ i \leq 5 $ such that $ {\zeta}^i \in f_0(w=t) $ or
                $ j \leq 3 $ such that $ y {\zeta}^j \in f_0(w=t), $ then
$(f_0(w=t),(0,0,0))$ is rational under the local coordinate system 
$(x,y,\zeta)$. 
So we assume the coefficient of $ {\zeta}^i $ is $0$ for all $ i \leq 5 $ and
             the coefficient of $ y {\zeta}^j $ is $0$ for all $ j \leq 3 $.
Furthermore we may assume $ \gamma \geq \delta $. ( Indeed, if 
$ \gamma < \delta $ we consider $ f_0(x,y,s,w) $ instead of $ f_0(x,y,z,t) $.) 
Thus $ b=b(t) $ is written as $ b(t)=b't^m, \; 
m := q/p = \gamma/\delta \geq 1, \; p, q \in {\bold N}, \; (p,q)=1 $. 
Since the coefficient of $ {\zeta}^i $ is
$ \displaystyle \frac{1}{i!}\frac{{\partial}^i f_0}{(\partial z)^i} (0,0,b,t) $
and the coefficient of $ y {\zeta}^j $ is $ \displaystyle 
\frac{1}{j!}\frac{{\partial}^{j+1} f_0}{\partial y (\partial z)^j} (0,0,b,t) $,
it follows that: \\
\hspace*{1.5cm} the coefficient of \ $ {\zeta}^i $ in $ f_0 (w=t) $ is $0$ \\
\hspace*{1cm} $ \Longleftrightarrow f_0 (0,0,z,t)=(z-b)^{i+1} \varphi(z,t) $, 
                for some $ \varphi(z,t) \in {\bold C} [z,t^m] $ \\
\hspace*{1cm} $ \Longleftrightarrow f_0 (0,0,z,w)=(z-b'w^m)^{i+1} \varphi(z,w)
              $, for some $ \varphi(z,w) \in {\bold C} [z,w^m] $,\\
and \\
\hspace*{1.5cm} the coefficient of \ $ y {\zeta}^j $ in $ f_0 (w=t) $ is $0$ \\
\hspace*{1cm} $ \displaystyle \Longleftrightarrow 
                \frac{\partial f_0}{\partial y}(0,0,z,t)=(z-b)^{j+1} \psi(z,t)
              $, for some $ \psi(z,t) \in {\bold C} [z,t^m] $ \\
\hspace*{1cm} $ \displaystyle \Longleftrightarrow 
                \frac{\partial f_0}{\partial y}(0,0,z,w)
                =(z-b'w^m)^{j+1} \psi(z,w) $, 
                for some $ \psi(z,w) \in {\bold C} [z,w^m] $.\\
The number $ m = \gamma/\delta $ is a integer by the assumptions. \\
( In fact, if $ p \ne 1 $ then $ f_0(0,0,z,w) $ and 
$ \displaystyle \frac{\partial f_0}{\partial y}(0,0,z,w) $ are written as \\
$ f_0(0,0,z,w)=(z-b'w^m)^6 \varphi(z,w)=(z^p-b''w^q)^6 \varphi'(z,w) , \,
b'' \in {\bold C}, \; \varphi'(z,w) \in {\bold C} [z,w] $,\\
$ \displaystyle \frac{\partial f_0}{\partial y}(0,0,z,w)=(z-b'w^m)^4 \psi(z,w)
=(z^p-b'''w^q)^4 \psi'(z,w) , $ \\
\hspace*{10cm} $ b''' \in {\bold C}, \; \psi'(z,w) \in {\bold C} [z,w] $,\\
respectively. Since $ 2 \leq p \leq q $, we have $ \displaystyle
\frac{1}{6p}+\frac{1}{6q} \leq \frac16 \quad \text{and \ }
\frac{1}{4p}+\frac{1}{4q} \leq \frac14 $. This is a contradiction to the
condition $ \gamma + \delta > 1/6 $.)

  Let $ z':=z-b'w^m $, \ $ \Lambda':=\{(k,l) \; ; {z'}^k w^l \in f_{\lambda}\} 
\bigcup\{\frac32(k,l) \; ; y {z'}^k w^l \in f_{\lambda}\} $, and let $\Gamma'$ 
be the union of the compact faces of the convex hull of $ \bigcup_{\nu \in 
\Lambda'} (\nu + {{\bold R}_{\geq 0}}^2) $ in $ {\bold R}^2 $. Then: 
\begin{claim} 
For any one-dimensional face $ \Delta' =\{ \gamma' Z + \delta' W = 1 , \; 
k'_1 \leq Z \leq k'_0, \; l'_0 \leq W \leq l'_1 \} $ of $\Gamma'$, 
the condition $ \gamma'+\delta'>1/6 $ is satisfied. 
\end{claim}

  A proof of this claim is found at the end of this paper. \\

  There exists $ k_1 < i \leq k_0 $ such that
$ i = \max \{ i' \in {\bold N} \; ; \; (z-b'w^m)^{i'} | f_0(0,0,z,w) \} $ 
or $ \frac23 k_1 < j \leq \frac23 k_0 $ such that
$ \displaystyle j = \max \left\{ j' \in {\bold N} \; ; \; (z-b'w^m)^{j'} 
\left| \frac{\partial f_0}{\partial y}(0,0,z,w) \right. \right\} $.
Choose $ \Delta' $ such that $ k'_0 \leq \min \{i, \; \frac32 j \} $, 
$ k'_1 < 6 $ and $ l'_0 < 6 $. Then $ k'_0 \leq k_0 $, \ $ l_0 \leq l'_0 $, 
and $ \displaystyle 1 \leq \frac{\gamma}{\delta} < \frac{\gamma'}{\delta'} $.  
Let $ f_0 $ be the initial part of $ f \in {\bold C} [x,y,z',w] $ with 
respect to the weight $(1/2, 1/3, \gamma', \delta')$.
If $ h:=f_0-x^2 \in {\bold C} [y,z',w] $  has a double factor, then 
Case (II-A). Now we assume $ h $ \ has no double factor.
If $ k'_0 < 6 $ then our assertion is concluded. 
So we assume $ k'_0 \geq 6 $. Repeating same argument as above, \\
\hspace*{1.5cm} $ z'':=z'-b''w^{m'}, \dotsb , 
                  z^{(n)}:=z^{(n-1)}-b^{(n)}w^{m^{(n-1)}}, \dotsb , $ \\
\hspace*{1.5cm} $ k_0 \geq k'_0 \geq k''_0 \geq \dotsb \geq k_0^{(n)} \geq 
                  \dotsb , $ \\
\hspace*{1.5cm} $ \displaystyle 1 \leq \frac{\gamma}{\delta} <
                  \frac{\gamma'}{\delta'} < \frac{\gamma''}{\delta''} < \dotsb
                  < \frac{{\gamma}^{(n)}}{{\delta}^{(n)}} < \dotsb , $ \\
\hspace*{1.5cm} $ \frac16 < {\gamma}^{(n)} + {\delta}^{(n)} $ \ for all 
                $ n \in {\bold N} $.\\ 
(See {\sc Figure} 2.) If there exists $ n \in {\bold N} $ such that 
$ k_0^{(n)} < 6 $, then this proof is completed. 
So consider the case $ k_0^{(n)} \geq 6 $ for all $ n \in {\bold N} $.  Then 
$ \displaystyle \frac{{\gamma}^{(n)}}{{\delta}^{(n)}} \in {\bold N} $ for all 
$ n \in {\bold N} $. If $ k_0^{(n)} > 6 $ for all $ n \in {\bold N} $ then  
$ \displaystyle \frac{{\gamma}^{(n)}}{{\delta}^{(n)}} 
            < - \frac{l_0^{(n)}-6}{k_0^{(n)}-6} $, 
and so this procedure must finish in finite times (See {\sc Figure} 3). 
Thus assume there exists $ n \in {\bold N} $ such that $ 6=k_0^{(n)}=
k_0^{(n+1)}=k_0^{(n+2)}=\dotsb . $ If this procedure continues infinitely then 
$ \mu(f_{\lambda}, 0) \gg 1 $ , a contradiction. This completes Case (I-A). \\
\begin{figure}[h]
\begin{center}
\setlength{\unitlength}{1mm}
\begin{picture}(90,80)(0,0)
  \put(70,10){\line(-1,1){40}}
  \put(60,20){\line(-2,3){24}}
  \put(50,35){\line(-1,4){10}}
  \put(70,5){\makebox(18,7){($k_0,l_0$)}}
  \put(12,45){\makebox(18,7){($k_1,l_1$)}}
  \put(60,20){\makebox(18,6){($k'_0,l'_0$)}}
  \put(17,55){\makebox(18,7){($k'_1,l'_1$)}}
  \put(50,33){\makebox(18,7){($k''_0,l''_0$)}}
  \put(22,75){\makebox(18,7){($k''_1,l''_1$)}}
  \put(45,25){$\Delta$}
  \put(50,55){$\Delta''$}
\end{picture}
\caption{ }
\end{center}
\end{figure}
\begin{figure}[h]
\begin{center}
\setlength{\unitlength}{1mm}
\begin{picture}(85,85)(-10,-10)
  \put(0,0){\vector(1,0){60}}
  \put(0,0){\vector(0,1){65}}
  \put(63,-2){$Z$}
  \put(-2,68){$W$}
  \put(-5,-5){0}
  \put(39,-5){6}
  \put(-5,39){6}
  \put(0,0){\dashbox(40,40)}
  \put(50,10){\line(-1,1){40}}
  \put(50,10){\line(-1,2){25}}
  \put(50,10){\line(-1,3){15}}
  \put(50,10){\makebox(20,7){($k_0^{(n)},l_0^{(n)}$)}}
\end{picture}
\caption{ }
\end{center}
\end{figure} \\
(I-B). \ $ C_j \ni P = (0,a(s),s,c(s)) \; ; \; s \ne 0 $.\\
Let $ \eta:=y-a, \; \omega:=w-c $, \; and $ f_0(z=s):=f_0(x, \eta+a, s,
\omega+c) $.
Then $ f_0(z=s)-x^2 \in {\bold C} [\eta,\omega] $ has no double factor.
If $ a \ne 0 $ or $ a=c=0 $ then $(f_0(x,y,s,w), \; (0,a,c))$ is rational, 
so $(f_0, \; (0,a,s,c))$ is rational similarly as in Case (I-A). Thus we 
assume $ a=0, \; c \ne 0 $. If $ \gamma \geq \delta $ then our argument can be 
reduced to the case $ a=0, \; b \ne 0 $ and $ \gamma \geq \delta $ of 
Case (I-A). So it is sufficient to consider the case $ \gamma < \delta $. 
Furthermore we may assume that the coefficient of  $ {\omega}^i $ in 
$ f_0(z=s) $ is $0$ for all $ i \leq 5 $ and the coefficient of  
$ y {\omega}^j $ in $ f_0(z=s) $ is $0$ for all $ j \leq 3 $. Thus $ c=c(s) $ 
is written as $ c(s)=c's^m, \; c' \in {\bold C}, \; m = \delta/\gamma > 1, \; 
m \in {\bold N} $. Therefore $ l_1 < 6 $ for No.11-14, 46-51 and 
$ l_1 \leq 6 $ for No.10, 83.
(Because $ z^{12} \in f $ for No.10,  $ z^{10} \in f $ for No.11,
          $ z^9 \in f $ for No.12,     $ z^8 \in f $ for No.13,
          $ z^7 \in f $ for No.14,     $ z^{11} \in f $ for No.46, 
          $ y z^7 \in f $ for No.47,   $ z^9 w \in f $ for No.48,
          $ z^8 w \in f $ for No.49,   $ z^7 w^2 \in f $ for No.50,
          $ z^7 w \in f $ for No.51,   $ z^{10} w \in f $ for No.83.)
Thus it is sufficient to consider the case $ (k_1,l_1)=(0,6) $ of No.10 and 
83. Let $ w':=w-c'z^m, \; \Lambda':=\{(k,l) \; ; z^k {w'}^l \in f_{\lambda}\} 
\bigcup \{\frac32(k,l) \; ; y z^k {w'}^l \in f_{\lambda}\} $, and let 
$\Gamma'$ be the union of the compact faces of the convex hull of 
$ \bigcup_{\nu \in \Lambda'} (\nu + {{\bold R}_{\geq 0}}^2) $ in ${\bold R}^2$.
Then for any one-dimensional face $ \Delta' =\{ \gamma' Z + \delta' W = 1 , \;
k'_1 \leq Z \leq k'_0, \; l'_0 \leq W \leq l'_1 \} $ of $\Gamma'$,
$ \gamma'+\delta'>1/6 $ and $ \displaystyle 1 < \frac{\delta}{\gamma} <
\frac{\delta'}{\gamma'} $ are satisfied. Repeating same argument as in Case 
(I-A), \\
\hspace*{1.5cm} $ w'':=w'-c''z^{m'}, \dotsb ,
                  w^{(n)}:=w^{(n-1)}-c^{(n)}z^{m^{(n-1)}}, \dotsb , $ \\
\hspace*{1.5cm} $ l_1 \geq l'_1 \geq l''_1 \geq \dotsb \geq l_1^{(n)} \geq
                  \dotsb , $ \\
\hspace*{1.5cm} $ \displaystyle 1 < \frac{\delta}{\gamma} <
                  \frac{\delta'}{\gamma'} < \frac{\delta''}{\gamma''} < \dotsb
                  < \frac{{\delta}^{(n)}}{{\gamma}^{(n)}} < \dotsb , $ \\
\hspace*{1.5cm} $ \frac16 < {\gamma}^{(n)} + {\delta}^{(n)} $ \ for all
                $ n \in {\bold N} $.\\
If there exists $ n \in {\bold N} $ such that $ l_1^{(n)} < 6 $ then the 
assertion holds true. If $ l_1^{(n)} = 6 $ for all \\ $ n \in {\bold N} $ 
then this procedure continues infinitely, therefore 
$ \mu(f_{\lambda}, 0) \gg 1 $, a contradiction. This completes Case (I).\\
\\
For convenience, we write $ z, \gamma, \delta, $ instead of 
$ z^{(n)}, {\gamma}^{(n)}, {\delta}^{(n)}, $ etc.\\
\\
Case II. \ $ h:=f_0-x^2 \in {\bold C} [y,z,w] $ \ has a double factor. \\
\hspace*{1cm} $ f_0 = x^2+y^3-\frac{27}{4}g(z,w)^2 y-\frac{27}{4}g(z,w)^3 
                    = x^2+(y+\frac32 g(z,w))^2 (y-3g(z,w)) $.\\
Let $ Y:= y+\frac32 g(z,w) $, then $f_0$ and $f_{\lambda}$ are written as:\\
\hspace*{1cm} $ f_0 = x^2+Y^3-\frac92 g(z,w) Y^2 $,\\
\hspace*{1cm} $ f_{\lambda} = x^2 + Y^3 - \frac92 g(z,w) Y^2 + \varphi(z,w) Y
                              + \psi(z,w) $, \\
for some $ \varphi(z,w), \; \psi(z,w) \in {\bold C} [z,w] $ with 
$ \varphi \not\equiv 0 $ or $ \psi \not\equiv 0 $, \ 
$ \varphi \equiv 0 $ or $ \deg_{\alpha}\varphi>2/3 $, \ 
$ \psi \equiv 0 $ or $ \deg_{\alpha}\psi>1 $.\\
\\
(II-A). \ When $ \gamma \geq \delta $, $g(z,w)$ can be classified into ten 
cases as below: \\
\hspace*{1cm} (1) \ \ $ g = c_{1 0} z + c_{0 L} w^L $, \\
\hspace*{1cm} (2) \ \ $ g = c_{1 1} z w + c_{0 (L+1)} w^{L+1} $, \\
\hspace*{1cm} (3) \ \ $ g = c_{2 0} z^2 + c_{0 L} w^L , \qquad 
                        3 \leq L $, \quad L is odd, \\
\hspace*{1cm} (4-a) \ $ g = c_{2 0} (z + \gamma_1 w^L)(z + \gamma_2 w^L) , 
                        \qquad \gamma_1 \ne \gamma_2 $,\\
\hspace*{1cm} (4-b) \ $ g = c_{2 0} (z + \gamma_1 w^L)^2 $, \\
\hspace*{1cm} (5) \ \ $ g = c_{2 1} z^2 w + c_{0 (L+1)} w^{L+1} , 
                        \qquad 3 \leq L $, \quad L is odd, \\
\hspace*{1cm} (6-a) \ $ g = c_{2 1} (z + \gamma_1 w^L)(z + \gamma_2 w^L) w,
                        \qquad \gamma_1 \ne \gamma_2 $,\\
\hspace*{1cm} (6-b) \ $ g = c_{2 1} (z + \gamma_1 w^L)^2 w $, \\
\hspace*{1cm} (7) \ \ $ g = c_{3 0} z^3 + c_{0 4} w^4 $, \\
\hspace*{1cm} (8) \ \ $ g = c_{3 0} z^3 + c_{1 3} z w^3 $, \\
\hspace*{1cm} (9) \ \ $ g = c_{3 0} z^3 + c_{0 5} w^5 $, \\
\hspace*{1cm} (10-a) \ $ g = c_{3 0} (z + \gamma_1 w)(z + \gamma_2 w)
                                        (z + \gamma_3 w) ,
                        \qquad \gamma_i \ne \gamma_j $ for $ i \ne j $,\\
\hspace*{1cm} (10-b) \ $ g = c_{3 0} (z + \gamma_1 w)^2 (z + \gamma_2 w) ,
                        \qquad \gamma_1 \ne \gamma_2 $,\\
\hspace*{1cm} (10-c) \ $ g = c_{3 0} (z + \gamma_1 w)^3 $. \\
After the local coordinate change $ z':=z+\gamma_1 w^L \; ( L=1 $ for (10))
around $0$, \\
\hspace*{1cm} (4-a) \ $ g = c_{2 0} {z'}^2 + c'_{1 L} z' w^L $, \\
\hspace*{1cm} (4-b) \ $ g = c_{2 0} {z'}^2 $, \\
\hspace*{1cm} (6-a) \ $ g = c_{2 1} {z'}^2 w + c'_{1 (L+1)} z' w^{L+1} $,\\
\hspace*{1cm} (6-b) \ $ g = c_{2 1} {z'}^2 w $, \\
\hspace*{1cm} (10-a) \ $ g = c_{3 0} {z'}^3 + c'_{2 1} {z'}^2 w 
                           + c'_{1 2} z' w^2 $, \\
\hspace*{1cm} (10-b) \ $ g = c_{3 0} {z'}^3 + c'_{2 1} {z'}^2 w $, \\
\hspace*{1cm} (10-c) \ $ g = c_{3 0} {z'}^3 $. \\
(See {\sc Figure} 4.)\\
\begin{figure}[h]
\begin{center}
\setlength{\unitlength}{1mm}
\begin{picture}(70,90)(-10,-10)
  \put(0,0){\vector(1,0){50}}
  \put(0,0){\vector(0,1){70}}
  \put(-5,-5){0}
  \put(9,-5){1}
  \put(19,-5){2}
  \put(29,-5){3}
  \put(39,-5){4}
  \put(-5,9){1}
  \put(-5,19){2}
  \put(-5,29){3}
  \put(-5,39){4}
  \put(-5,49){5}
  \put(-5,59){6}
  \put(53,-2){$Z$}
  \put(-2,73){$W$}
  \put(20,20){\circle*{1}}
  \put(20,20){\makebox(12,6){(2,2)}}
  \put(10,0){\line(-1,3){10}}
  \put(10,10){\line(-1,5){10}}
  \put(20,0){\line(-1,1){10}}
  \put(20,10){\line(-1,3){10}}
  \put(30,0){\line(-3,4){30}}
  \put(30,0){\line(-2,3){20}}
  \put(30,0){\line(-3,5){30}}
  \multiput(30,0)(-1,1){20}{\circle*{0.2}}
  \put(20,0){\line(-2,3){20}}
  \put(10,0){\circle*{1}}  
  \put(20,0){\circle*{1}}
  \put(30,0){\circle*{1}}
  \put(10,10){\circle*{1}}
  \put(20,10){\circle*{1}}
  \put(0,30){\circle*{1}}
  \put(10,30){\circle*{1}}
  \put(0,40){\circle*{1}}
  \put(10,40){\circle*{1}}
  \put(0,50){\circle*{1}}
  \put(0,60){\circle*{1}}
  \put(10,20){\circle{1}}
\end{picture}
\caption{ }
\end{center}
\end{figure} 
Let $ \Lambda':=\{(k,l) \; ; {z'}^k w^l \in f_{\lambda}\}
\bigcup\{\frac32(k,l) \; ; Y {z'}^k w^l \in f_{\lambda}\}
\bigcup\{3(k,l) \; ; Y^2 {z'}^k w^l \in f_{\lambda}\} $, and let $\Gamma'$ be
the union of the compact faces of the convex hull of $ \bigcup_{\nu \in
\Lambda'} (\nu + {{\bold R}_{\geq 0}}^2) $ in $ {\bold R}^2 $. 
Then for any one-dimensional face $ \Delta' =\{ \gamma' Z + \delta' W = 1 \; ,
\; k'_1 \leq Z \leq k'_0, \; l'_0 \leq W \leq l'_1 \} $ of $\Gamma'$, the 
condition $ \gamma'+\delta'>1/6 $ is satisfied. (In fact, No.11, 13, 14 and 
49-51 can not become Case (10). About Case (10-c) which comes from No.10, 12, 
46-48 or 83, the condition $ \gamma'+\delta'>1/6 $ is satisfied as in 
Claim 3.5.) \\
\\
First we consider Case (1), (2), (3), (4-a), (5), (6-a), (7), (8), (9) and 
(10-a). \\
We replace the weight $ \alpha=(1/2,1/3,\gamma,\delta) $ with
$ \alpha'=(1/2,\beta',\gamma',\delta') $ which satisfies the \\ conditions \\
\hspace*{1cm} $ \frac13 = \beta < \beta' < \frac12 
                                         < \beta' + \gamma' + \delta', \quad
                \gamma'/\delta' = \gamma/\delta, $ \\
\hspace*{1cm} $ \deg_{\alpha'}(g(z,w)Y^2) = 1, $ \; 
              $ \deg_{\alpha'}(\varphi(z,w)Y), \;  
                \deg_{\alpha'}\psi(z,w) \geq 1 $, \; and \\
\hspace*{1cm} $ \deg_{\alpha'}(\varphi(z,w)Y) $ or 
              $ \deg_{\alpha'}\psi(z,w) = 1 $,\\
and let $ f_0 $ be the initial part of $ f_{\lambda} $.
(There exists such weight $\alpha'$ because \\
\hspace*{1cm} $   \beta' + \gamma' + \delta' 
                = \beta' + (\gamma+\delta)(1-2\beta')/(1-2\beta) > \frac12 $ 
              and \\
\hspace*{1cm} $ \deg_{\alpha'}(g(z,w)Y^2) 
              = 2\beta'+\deg_{\alpha}g(z,w)(1-2\beta')/(1-2\beta) = 1 $ \\
for all $ \beta, \gamma, \delta, \beta', \gamma', \delta' $ satisfying \\
\hspace*{1cm} $ \beta, \; \beta' < \frac12 < \beta+\gamma+\delta $ and 
              $ \gamma'/\gamma=\delta'/\delta=(1-2\beta')/(1-2\beta) $.) \\
Then \\
\hspace*{1cm} $ f_0 = x^2 - \frac92 g(z,w) Y^2
                      + \varphi_0 (z,w) Y + \psi_0 (z,w) $,\\
\hspace*{1cm} $ f_{\lambda} = x^2 + Y^3 - \frac92 g(z,w) Y^2
                      + \varphi (z,w) Y + \psi (z,w) $,\\
where $ \varphi_0 \; (\text{resp. } \psi_0) $ is either 0 or the initial part 
of $ \varphi \; (\text{resp. }\psi) $. \\
\\
(II-A-I). \ $ h:=f_0-x^2 \in {\bold C} [Y,z,w] $ \ has no double factor. \\
Let $ C_j \subset \operatorname{Sing}(f_0) $ be an irreducible curve with $ 0 
\in C_j $, and $ P = (0,a,b,c) \ne 0 $ an arbitrary point on $C_j$. \\
\\
(II-A-I-i). \ $ C_j \ni P = (0,a(t),b(t),t) \; ; \; t \ne 0 $.\\
Then, $ a(t)=a't^{\beta'/\delta'}, \; b(t)=b't^{\gamma'/\delta'} $ for some 
$ a', \; b' \in {\bold C} $.\\
Let $ \eta:=Y-a, \; \zeta:=z-b $, \ and $ f_0(w=t):=f_0(x, \eta+a, \zeta+b, t) 
$, then \\
\hspace*{1cm} $ f_0(w=t) = x^2 - \frac92 g(\zeta+b,t) (\eta+a)^2
                           + \varphi_0 (\zeta+b,t) (\eta+a) 
                           + \psi_0 (\zeta+b,t) $ \\
\hspace*{2.8cm} $ = x^2 - \frac92 g(\zeta+b,t) \eta^2
                    + (-9a g(\zeta+b,t) + \varphi_0 (\zeta+b,t)) \eta $ \\
\hspace*{3.2cm}   $ - \frac92 a^2 g(\zeta+b,t) + a \varphi_0 (\zeta+b,t) 
                    + \psi_0 (\zeta+b,t) $.\\
Then both $ g(\zeta+b,t) \in {\bold C} [\zeta] $ and 
$ f_0(w=t)-x^2 \in {\bold C} [\eta,\zeta] $ have no double factor.\\
Therefore $ (f_0(w=t), \; (0,0,0)) $ is an isolated singularity under the 
coordinates $ (x, \eta, \zeta) $. \\
If $ g(b,t) \ne 0 $ then $ {\eta}^2 \in f_0(w=t) $. 
Otherwise, $ \zeta \in g(\zeta+b,t) $, so $ {\eta}^2 \zeta \in f_0(w=t) $. \\
Hence $ (f_0(w=t), \; (0,0,0)) $ is rational under the coordinates 
$ (x, \eta, \zeta) $ by Lemma 3.2. \\
\\
(II-A-I-ii). \ $ C_j \ni P = (0,a(s),s,0) \; ; \; s \ne 0 $.\\
Let $ \eta:=Y-a $, then \\
\hspace*{1cm} $ f_0(z=s) = x^2 - \frac92 g(s,w) (\eta+a)^2 
                           + \varphi_0 (s,w) (\eta+a) + \psi_0 (s,w) $ \\
\hspace*{2.8cm} $ = x^2 - \frac92 g(s,w) {\eta}^2 
                    + ( -9 a g(s,w) + \varphi_0 (s,w)) \eta $ \\
\hspace*{3.3cm} $ - \frac92 a^2 g(s,w) + a \varphi_0 (s,w) + \psi_0 (s,w). $ \\
$ (f_0(z=s), \; (0,0,0)) $ is an isolated singularity under the coordinates
$ (x, \eta, w) $ because $ h:=f_0(z=s)-x^2 \in {\bold C}[\eta,w] $ has no 
double factor. We have 
$ {\eta}^2 \in f_0(z=s) $ for (1), (3), (4-a), (7), (8), (9), (10-a), and 
$ {\eta}^2 w \in f_0(z=s) $ for (2), (5), (6-a), 
since $ s \ne 0 $. Therefore $ (f_0(z=s), \; (0,0,0)) $ is rational. \\
\\
(II-A-I-iii). \ $ C_j \ni P = (0,a,0,0) \; ; \; a \ne 0 $.\\
Since $ f_0(Y=a) = x^2 - \frac92g(z,w) a^2 + \varphi_0(z,w) a + \psi_0(z,w) $, 
it follows that \\ 
$ (f_0(Y=a), \; (0,0,0)) $ is rational under the coordinates $ (x,z,w) $ 
by Lemma 3.2.\\
\\
(II-A-II). \ $ h:=f_0-x^2 \in {\bold C}[Y,z,w] $ has a double factor. \\
\hspace*{1cm} $ f_{\lambda} = x^2 + Y^3 - \frac92 g(z,w) Y^2
                              + \varphi (z,w) Y + \psi (z,w) $,\\
\hspace*{1cm} $ f_0 = x^2 - \frac92 g(z,w) Y^2
                      + \varphi_0 (z,w) Y + \psi_0 (z,w) $ \\
\hspace*{1.5cm}   $ = x^2 - \frac92 g(z,w) (Y+\phi(z,w))^2 $.\\
After the coordinate change $ Y':=Y+\phi(z,w) $, \\
\hspace*{1cm} $ f_{\lambda} = x^2 + {Y'}^3 + (-\frac92 g - 3\phi)(z,w){Y'}^2
                              + \varphi'(z,w) Y'+ \psi'(z,w) $ \\ 
for some $ \varphi', \; \psi' \in {\bold C}[z,w] $ with 
$ \deg_{\alpha'} \varphi < \deg_{\alpha'} \varphi', \; 
  \deg_{\alpha'} \psi < \deg_{\alpha'} \psi' $. 
We replace the weight $ \alpha'=(1/2,\beta',\gamma',\delta') $ with $ \alpha''
=(1/2,\beta'',\gamma'',\delta'') $ which satisfies the conditions \\
\hspace*{1cm} $ \frac13 < \beta'' < \frac12 < \beta''+\gamma''+\delta'', \quad
                \gamma''/\delta'' = \gamma'/\delta', $ \\
\hspace*{1cm} $ \deg_{\alpha''}(g(z,w){Y'}^2) = 1, \quad 
                \deg_{\alpha''}(\varphi'(z,w) Y'), \; 
                \deg_{\alpha''}\psi'(z,w) \geq 1 $, \; and \\
\hspace*{1cm} $ \deg_{\alpha''}(\varphi'(z,w) Y') \; \text{or} \; 
                \deg_{\alpha''}\psi'(z,w) = 1 $,\\
and let $ f_0 $ be the initial part of $ f_{\lambda} $. Then \\
\hspace*{1cm} $ f_0 = x^2 - \frac92 g(z,w) {Y'}^2
                      + \varphi'_0 (z,w) Y' + \psi'_0 (z,w) $, \\
\hspace*{1cm} $ \beta' < \beta'', \quad 
                \gamma' > \gamma'', \quad
                \delta' > \delta'', $ \\
\hspace*{1cm} $ \deg_{\alpha'}g(z,w) < \deg_{\alpha'}\phi(z,w), $ \\
\hspace*{1cm} $ \deg_{\alpha'}\varphi_0(z,w) < \deg_{\alpha'}\varphi'_0(z,w) $ 
              or $ \varphi'_0(z,w) \equiv 0, $ \\
\hspace*{1cm} $ \deg_{\alpha'} \psi_0(z,w) < \deg_{\alpha'} \psi'_0(z,w) $ or 
              $ \psi'_0(z,w) \equiv 0. $ \\
If $ h:=f_0-x^2 \in {\bold C}[Y',z,w] $ has no double factor, then Case 
(II-A-I). \\
If $h$ has a double factor, namely, \\
\hspace*{1.5cm} $ f_0 = x^2 - \frac92 g(z,w) (Y'+\phi'(z,w))^2 $,\\
then after the coordinate change $ Y'':=Y'+\phi'(z,w) $, \\
\hspace*{1cm} $ f_{\lambda} = x^2 + {Y''}^3 
                              + (-\frac92 g - 3(\phi+\phi'))(z,w){Y''}^2
                              + \varphi''(z,w) Y''+ \psi''(z,w) $ \\
\hspace*{1.5cm} $ =: x^2 + {Y''}^3 + (-\frac92 g - 3\phi)(z,w){Y''}^2
                     + \varphi''(z,w) Y''+ \psi''(z,w) $ \\ 
for some $ \varphi'', \; \psi'' \in {\bold C}[z,w] $ with
$ \deg_{\alpha''} \varphi' < \deg_{\alpha''} \varphi'', \;
  \deg_{\alpha''} \psi' < \deg_{\alpha''} \psi'' $.\\
If this procedure continues infinitely, then \\
\hspace*{1cm} $ f_{\lambda} = x^2 + {Y^{(n)}}^3 
                              + (-\frac92 g - 3\phi)(z,w){Y^{(n)}}^2
                              + \varphi^{(n)}(z,w)Y^{(n)}+\psi^{(n)}(z,w) $,\\
\hspace*{1cm} $ \frac13 = \beta < \beta' < \beta'' < \dotsb < \beta^{(n)} 
                < \dotsb < \frac12, $ \\
\hspace*{1cm} $ \gamma > \gamma' > \gamma'' > \dotsb > \gamma^{(n)} 
                > \dotsb, $ \\
\hspace*{1cm} $ \delta > \delta' > \delta'' > \dotsb > \delta^{(n)} 
                > \dotsb, $ \\
\hspace*{1cm} $ \gamma/\delta = \gamma'/\delta' = \gamma''/\delta'' = \dotsb = 
                \gamma^{(n)}/\delta^{(n)} = \dotsb, $ \\
\hspace*{1cm} $ \frac12 < \beta^{(n)} + \gamma^{(n)} + \delta^{(n)} \quad 
                \text{for all } n \in {\bold Z}_{\geq 0}, $ \\
\hspace*{1cm} $ \deg_{\alpha} g(z,w) < \deg_{\alpha} \phi(z,w), $ \\
\hspace*{1cm} $ \deg_{\alpha} \varphi_0(z,w) < \deg_{\alpha} \varphi'_0(z,w) 
                < \deg_{\alpha} \varphi''_0(z,w) < \dotsb 
                < \deg_{\alpha} \varphi_0^{(n)}(z,w) < \dotsb, $ \\
\hspace*{1cm} $ \deg_{\alpha} \psi_0(z,w) < \deg_{\alpha} \psi'_0(z,w) 
                < \deg_{\alpha} \psi''_0(z,w) < \dotsb 
                < \deg_{\alpha} \psi_0^{(n)}(z,w) < \dotsb, $ \\
so $ \displaystyle 
   \mu(f_{\lambda},0) \geq 
   \dim_{\bold C} {\bold C} \{ w \} \left/
   \left( \varphi^{(n)}(0,w), \;
   \frac{\partial \psi^{(n)}}{\partial z}(0,w), \;
   \frac{\partial \psi^{(n)}}{\partial w}(0,w) \right) \right. \gg 1 $,\\
a contradiction. \\
\\
Next we consider Case (4-b), (6-b), (10-b), (10-c).
Choose a $\Delta'$ such that: \\
\hspace*{1cm} $ (k'_0, l'_0) $ = (6, 0) for (4-b), \ 
                                 (6, 3) for (6-b) and (10-b), \\
\hspace*{2.8cm} $ k'_1<6 $ and $ l'_0<6 $ for (10-c), \\
respectively. 
Then $ \displaystyle 1 \leq \frac{\gamma}{\delta} < \frac{\gamma'}{\delta'} $. 
Let $f_0$ be the initial part of $ f_{\lambda} \in {\bold C} [x,Y,z',w] $
with respect to the weight $ \alpha':=(1/2, 1/3, \gamma', \delta') $, and 
$ m':=\gamma'/\delta' $.\\
For the case ${z'}^3 Y^2$ is not contained in $f_0$ of (10-c), the situation is
similar to (I-A), but $ h:=f_0-x^2 \in {\bold C}[Y,z^{(n)},w] $ can not have a 
double factor for $ n \in {\bold N} $. So we consider other cases, i.e., 
\begin{align*}
\hspace*{1cm} f_0 &= x^2 + Y^3 - \frac92c_{kl}^{(\prime)}{z'}^k w^l Y^2 
                     + \varphi'_0(z',w)Y + \psi'_0(z',w) \\
                  &=
 \begin{cases}
  x^2 + Y^3 - \frac92c_{20}{z'}^2 Y^2 + \varphi'_0(z',w)Y + \psi'_0(z',w), 
  \qquad & \text{(4-b)} \\
  x^2 + Y^3 - \frac92c_{21}{z'}^2 w Y^2 + \varphi'_0(z',w)Y + \psi'_0(z',w), 
  \qquad & \text{(6-b)} \\
  x^2 + Y^3 - \frac92c'_{21}{z'}^2 w Y^2 + \varphi'_0(z',w)Y + \psi'_0(z',w), 
  \qquad & \text{(10-b)} \\
  x^2 + Y^3 - \frac92c_{30}{z'}^3 Y^2 + \varphi'_0(z',w)Y + \psi'_0(z',w), 
  \qquad & \text{(10-c)}
 \end{cases}
\end{align*} 
where $ \deg\varphi'_0(z',t)<2k, \; \deg\psi'_0(z',t)<3k \quad 
(0 \ne t \in {\bold C}). $ \\
\\
(II-A-I$'$). \ $ h:=f_0-x^2 \in {\bold C} [Y,z',w] $ \ has no double factor. \\
(II-A-I$'$-i). \ $ \text{Sing}(f_0) \supset C_j \ni P=(0,a(t),b(t),t) \; ; \; 
                  t \ne 0 $.\\
Let $ \eta:=Y-a, \; \zeta:=z'-b $ and $ f_0(w=t):=f_0(x, \eta+a, \zeta+b, t) $,
then $ f_0(w=t)-x^2 $ has no double factor. \\
(Indeed, assume $ f_0(w=t)-x^2=(\eta+A(\zeta))^2(\eta+B(\zeta)) $ for some 
$ A, \; B \in {\bold C}[\zeta] $, then \\
\hspace*{1cm} $ f_0(w=t)-x^2=(Y-a(t)+A(z'-b(t)))^2(Y-a(t)+B(z'-b(t))) $,\\
\hspace*{1cm} $ f_0-x^2=(Y-a(w)+A(z'-b(w)))^2(Y-a(w)+B(z'-b(w))) $ \\
\hspace*{2.4cm} $ =(Y+G(z',w))^2(Y+H(z',w)) $ \\
\hspace*{2.4cm} $ (G:=-a(w)+A(z'-b(w)), \; H:=-a(w)+B(z'-b(w)) 
                  \in {\bold C}[z',w^{m'}].) $ \\
\hspace*{2.4cm} $ =Y^3+(2G+H)Y^2+G(G+2H)Y+G^2H $ \\
\hspace*{2.4cm} $ =Y^3-\frac92g'Y^2-G(9g'+3G)Y-G^2(\frac92g'+2G) $,\\
where $ g'(z',w):=c_{k l}^{(\prime)}{z'}^k w^l $. 
It follows that $ \deg G(z',t) < k $ from $ \deg \varphi'_0(z',t) < 2k, \; 
\deg \psi'_0(z',t) < 3k $, and that $ G \in {\bold C}[z',w] $ from $ g', \; 
G(9g'+3G), \; G^2(\frac92g'+2G) \in {\bold C}[z',w], \; 
G \in {\bold C}[z',w^{m'}] $.)\\
If $ a=b=0 $ then \\
\hspace*{2cm} $ f_0(w=t) = x^2 + Y^3 - \frac92c_{k l}^{(\prime)}{z'}^k t^l Y^2
                + \varphi'_0(z',t)Y + \psi'_0(z',t), $ \\
with $ \text{ord}\varphi'_0(z',t) < 4 \leq 2k $ or 
     $ \text{ord}\psi'_0(z',t) < 6 \leq 3k $.\\
If $ a=0, \; b \ne 0 $, then \\
\hspace*{1cm} $ f_0(w=t) = x^2 + Y^3 - \frac92c_{k l}^{(\prime)}
                           ({\zeta}^k+kb{\zeta}^{k-1} 
                           +\dotsb+kb^{k-1}\zeta+b^k)t^l Y^2 $ \\
\hspace*{3.3cm} $ + \varphi'_0(\zeta+b,t)Y + \psi'_0(\zeta+b,t) \ni Y^2. $ \\
If $ a \ne 0, \; b=0 $, then \\
\hspace*{1cm} $ f_0(w=t) = x^2 + {\eta}^3 + 3a{\eta}^2 
                + (3a^2 + \varphi'_0(z',t))\eta 
                - \frac92c_{k l}^{(\prime)}{z'}^k t^l (\eta + a)^2 + \dotsb 
                \ni {\eta}^2. $ \\
If $ a \ne 0, \; b \ne 0 $, then $ {\eta}^2 \in f_0(w=t) $.

(Indeed, the coefficient of ${\eta}^2$ in $f_0(w=t)$ is $ \displaystyle 
\frac12 \frac{\partial^2 f_0}{(\partial Y)^2}(0,a,b,t) $, so it follows that:\\
\hspace*{2.5cm} the coefficient of ${\eta}^2$ in $f_0(w=t)$ is 0 \\
\hspace*{2cm} $ \Longleftrightarrow f_0(0,Y,b(t),t)=(Y-a(t))^3 $ \\
\hspace*{2cm} $ \Longleftrightarrow f_0(0,Y,b(w),w)=(Y-a(w))^3 $ \\
\hspace*{2cm} $ \Longleftrightarrow 
                f_0(0,Y,z',w)=(Y-\frac32c_{k l}^{(\prime)}{z'}^k w^l)^3 $.)\\
Hence $(f_0(w=t), 0)$ is rational. \\
\\
(II-A-I$'$-ii). \ $ \text{Sing}(f_0) \supset C_j \ni P=(0,a(s),s,0) \; ; \; 
                    s \ne 0 $.\\
Then $ a=0 $ because 
$ \displaystyle 0=f_0(0,a,s,0)=\frac{\partial f_0}{\partial Y}(0,a,s,0)
                 =\varphi'_0(s,0)=\psi'_0(s,0) $.\\
Furthermore, \\
\hspace*{1cm} $ f_0(z'=s) - x^2 :=f_0(x, Y, s, w) - x^2 
                                 = Y^3 - \frac92c_{k l}^{(\prime)}s^k w^l Y^2
                                   + \varphi'_0(s,w)Y + \psi'_0(s,w) $ \\
has no double factor. So it follows that $(f_0(z'=s),0)$ is rational from \\
\hspace*{1cm} $ Y^2 \in f_0(z'=s) $ for (4-b) and (10-c), \\
\hspace*{1cm} $ Y^2 w \in f_0(z'=s) $ and $ w^i, \; Y w^j $ are not contained 
in $f_0(z'=s)$ for $ i \leq 3, \; j \leq 2 $ \\
\hspace*{1cm} for (6-b) and (10-b).\\
\\
(II-A-II$'$). \ $ h:=f_0-x^2 \in {\bold C} [Y,z',w] $ \ has a double factor. \\
$f_0$ and $f_{\lambda}$ are written as: \\
\hspace*{1cm} $ f_0=x^2+(Y+G(z',w))^2 (Y+(-\frac92g'-2G)(z',w)) $ \\
\hspace*{1.5cm} $ =:x^2+{Y'}^3-(\frac92g'+3G)(z',w){Y'}^2, \qquad 
                  (Y':=Y+G(z',w).) $ \\
\hspace*{1cm} $ f_{\lambda}=:
  \begin{cases}
    x^2+{Y'}^3-\frac92g''(z',w){Y'}^2+\varphi''(z',w)Y'+\psi''(z',w), \qquad \;
    \text{(4-b, 6-b, 10-c)} \\
    x^2+{Y'}^3-\frac92(c_{3 0}{z'}^3+g''(z',w)){Y'}^2+\varphi''(z',w)Y'
    +\psi''(z',w), \qquad \text{(10-b)} 
  \end{cases} $ \\
\\
for some $ g'':=g'+\frac23G, \; \varphi'', \; \psi'' \in {\bold C}[z',w] $ with
$ \deg G(z',t) < k \quad (0 \ne t \in {\bold C}), \\ 
  \deg_{\alpha'}\varphi''>2/3=2(k\gamma'+l\delta'), \; 
  \deg_{\alpha'}\psi''>1 $.\\
Then $g''$ is one of the following:\\
\\
$ \hspace*{1cm} \left. \begin{gathered}
   \text{(3) \ } g'' = c_{2 0} {z'}^2 + c'_{0 L'} w^{L'}, \qquad 
                   3 \leq L', \quad L' \; \text{is odd}, \hspace{1.2cm} \\
   \text{(4-a) \ } g''= c_{2 0}(z'+ \gamma'_1 w^{L'})(z'+ \gamma'_2 w^{L'}), 
                    \quad 2 \leq L', \; \gamma'_1 \ne \gamma'_2, \\
   \text{(4-b) \ } g'' = c_{2 0} (z' + \gamma'_1 w^{L'})^2, \qquad 2 \leq L',
                     \hspace{3.2cm}
 \end{gathered} \right\} \quad \text{(become from 4-b)} \\
\\
 \hspace*{1cm} \left. \begin{gathered}
   \text{(5) \ } g'' = c_{2 1}^{(\prime)}{z'}^2 w + c'_{0 (L'+1)}w^{L'+1}, 
                          \qquad 3 \leq L', \quad L' \; \text{is odd}, \; \\
   \text{(6-a) \ } g'' = c_{2 1}^{(\prime)} (z' + \gamma'_1 w^{L'})
                                              (z' + \gamma'_2 w^{L'}) w,
                     \quad 2 \leq L', \; \gamma'_1 \ne \gamma'_2, \\
   \text{(6-b) \ } g'' = c_{2 1}^{(\prime)} (z' + \gamma'_1 w^{L'})^2 w, 
                      \qquad 2 \leq L', \hspace{3.2cm}
 \end{gathered} \right\} \quad  
 \begin{gathered}
   \text{(become from} \\
   \text{6-b or 10-b)} 
 \end{gathered} \\
\\
 \hspace*{1cm} \left. \begin{gathered}
   \text{(7) \ } g'' = c_{3 0} {z'}^3 + c'_{0 4} w^4, \quad \\
   \text{(8) \ } g'' = c_{3 0} {z'}^3 + c'_{1 3} z' w^3, \\
   \text{(9) \ } g'' = c_{3 0} {z'}^3 + c'_{0 5} w^5, \quad
 \end{gathered} \right\} \quad \text{(become from 10-c).} $ \\
\\
For 3, 4-a, 5, 6-a, 7, 8 and 9, we replace the weight 
$ \alpha'=(1/2, 1/3, \gamma', \delta') $ with 
$ \alpha''=(1/2, \beta'', \gamma'', \delta'') $ 
which satisfies the conditions \\
\hspace*{1cm} $ \frac13 < \beta'' < \frac12 < \beta''+\gamma''+\delta'', \quad
                \gamma''/\delta'' = \gamma'/\delta', $ \\
\hspace*{1cm} $ \deg_{\alpha''}(g''(z',w){Y'}^2) = 1, \quad 
                \deg_{\alpha''}(\varphi''(z',w) Y'), \; 
                \deg_{\alpha''}\psi''(z',w) \geq 1 $, \; and \\
\hspace*{1cm} $ \deg_{\alpha''}(\varphi''(z',w) Y') \; \text{or} \; 
                \deg_{\alpha''}\psi''(z',w) = 1 $,\\
and let $ f_0 $ be the initial part of $ f_{\lambda} $.\\ 
If $ h:=f_0-x^2 \in {\bold C}[Y',z',w] $ has no double factor then Case 
(II-A-I). Otherwise, Case (II-A-II), and this procedure must finish in finite 
times from the assumption. \\
For 4-b, 6-b, after the coordinate change $z'':=z'+\gamma'_1 w^{L'}$, 
we define $\Lambda''$, $\Gamma''$ under the coordinates $(x, Y', z'', w)$ 
similarly as before and choose a $\Delta''$ corresponding to a weight
$ \alpha''=(1/2, 1/3, \gamma'', \delta'') $ which satisfies $ (k''_0, l''_0)
=(6, 0) $ for (4-b), \ (6, 3) for (6-b). Then $ \displaystyle \frac16 < 
\gamma''+\delta''$ and $ \displaystyle 1 \leq \frac{\gamma}{\delta} < 
\frac{\gamma'}{\delta'} < \frac{\gamma''}{\delta''} $. Let $f_0$ be the initial
part of $f_{\lambda}$. Then, \\
\\
\hspace*{1cm} $ f_0 = 
  \begin{cases}
    x^2+{Y'}^3-\frac92 c_{2 0}{z''}^2{Y'}^2+\varphi'''_0(z'',w)Y'
    +\psi'''_0(z'',w), \qquad & \text{(4-b)} \\
    x^2+{Y'}^3-\frac92 c_{2 1}^{(\prime)}{z''}^2 w{Y'}^2+\varphi'''_0(z'',w)Y'
    +\psi'''_0(z'',w), \qquad & \text{(6-b)} 
  \end{cases} $ \\
\\
where $ \deg\varphi'''_0(z'',t) < 4, \; \deg\psi'''_0(z'',t) < 6 \quad 
(0 \ne t \in {\bold C}) $. \\  
If $ h:=f_0-x^2 \in {\bold C}[Y',z'',w] $ has no double factor 
then Case (II-A-I$'$). Otherwise, Case (II-A-II$'$), and this procedure must 
finish in finite times from the assumption. \\
\\
(II-B). \ When $ \gamma < \delta $, $ g(z,w) $ can be classified as below: \\
\hspace*{1cm} (1$'$) \ \ $ g = c_{0 1} w + c_{K 0} z^K, \qquad 
                         2 \leq K $,\\
\hspace*{1cm} (2$'$) \ \ $ g = c_{1 1} z w + c_{(K+1) 0} z^{K+1}, \qquad
                         2 \leq K $, \\
\hspace*{1cm} (3$'$) \ \ $ g = c_{0 2} w^2 + c_{K 0} z^K, \qquad 
                         3 \leq K $, \  K is odd, \\
\hspace*{1cm} (4-a$'$) \ $ g = c_{0 2} (w + \gamma_1 z^K)(w + \gamma_2 z^K), 
                         \quad 2 \leq K, \ \gamma_1 \ne \gamma_2 $, \
                         ( for No.10, 83 only ), \\
\hspace*{1cm} (4-b$'$) \ $ g = c_{0 2} (w + \gamma_1 z^K)^2, \qquad 
                         2 \leq K $, \qquad ( for No.10, 83 only ), \\
\hspace*{1cm} (5$'$) \ \ $ g = c_{1 2} z w^2 + c_{(K+1) 0} z^{K+1}, \quad
                         3 \leq K $, \  K is odd, \qquad 
                         ( for No.10 only ),\\
\hspace*{1cm} (7$'$) \ \ $ g = c_{0 3} w^3 + c_{4 0} z^4 $, \qquad
                         ( for No.10 only ). \\
For Cases (1$'$), (2$'$), (3$'$), (4-a$'$), (5$'$) and (7$'$), the situation 
is similar to (II-A).
For Case (4-b$'$) of No.10 and 83, we define $ \Lambda', \Gamma', \Delta', f_0 
$ and $h$ similarly as in (II-A), and repeat the same argument as in (II-A). 
Namely, if $ h:=f_0-x^2 \in {\bold C}[Y,z,w'] $ has no double factor then 
(II-B-I$'$), and if $h$ has a double factor then (II-B-II$'$). 
Then the condition $ \beta^{(n)}+\gamma^{(n)}+\delta^{(n)}>\frac12 $ must be 
always satisfied for No.10, because \\ $ \mu(f_{\lambda},0)<\mu(f,0)=242 $. 
For the case of $(k_1,l_1)=(k'_1,l'_1)=(0,6)$ of No.83, if it becomes 
(II-B-II$'$), then $ z^k \in f_{\lambda} $ for some $ k \geq 15 $. 
And so, if $ \beta^{(n)}+ \gamma^{(n)}+\delta^{(n)} \leq 1/2 $ for some 
$ n \in {\bold N} $ then $ \mu(f_{\lambda},0) \geq (3-1)(10-1)(15-1)=252 $, 
which is a contradiction to the condition $ \mu(f_{\lambda},0)<\mu(f,0)=245 $.

And this argument must finish finitely since $ \mu(f_{\lambda},0)<\mu(f,0) $.\\

For arrangement, let us illustrate the above argument as follows:\\
\\
{\small $
 \text{I} 
 \begin{cases}
  \text{I-A}
  \begin{cases}
   a \ne 0 \quad \; \; : \text{O.K.} \dotsb \; (*) \\
   a=b=0 : \text{O.K.} \dotsb \; (**) \\
   a=0, \; b \ne 0 
   \begin{cases}
    \zeta^i \; (\exists \; i \leq 5) \; \text{or} \; 
    y\zeta^j \; (\exists \; j \leq 3) \in f_0(w=t) : 
    \text{O.K.} \dotsb \; (\text{***}) \\
    \text{otherwise} : z':=z-b'w^m \; \longrightarrow 
    \begin{cases}
     \text{I-A} \\  
     \text{II-A}
    \end{cases}
   \end{cases}
  \end{cases}\\
  \\
  \text{I-B} 
  \begin{cases}
   a \ne 0 \quad \; \; : \text{O.K.} \dotsb \; (*') \\
   a=c=0 : \text{O.K.} \dotsb \; (**') \\
   a=0, \; c \ne 0 
   \begin{cases}
    \omega^i \; (\exists \; i \leq 5) \; \text{or} \; 
    y\omega^j \; (\exists \; j \leq 3) \in f_0(z=s) : 
    \text{O.K.} \dotsb \; (\text{***}') \\
    \text{otherwise} : w':=w-c'z^m \; \longrightarrow 
    \begin{cases}
     \text{I-B} \\
     \text{II-B}
    \end{cases}
   \end{cases}
  \end{cases}   
 \end{cases} \\
\\
 \text{II}
 \begin{cases}
  \text{II-A} 
  \begin{cases}
   \begin{gathered}
    \text{1, 2, 3, 4-a, 5, } \\
    \text{6-a, 7, 8, 9, 10-a }  
   \end{gathered} 
   \begin{cases}
    \text{II-A-I : O.K.} \dotsb \; (*'') \\
    \text{II-A-II} : Y':=Y+\phi \; \to 
    \begin{cases}
     \text{II-A-I : O.K.} \dotsb \; (*'') \\
     \text{II-A-II} : Y'':=Y'+\phi' \; \to \; \dotsb 
    \end{cases}
   \end{cases} \\
   \\
   \begin{gathered}
    \text{4-b} \\
    \text{6-b} \\
    \text{10-b} 
   \end{gathered} : z':=z+\gamma_1 w^L \; \to 
   \begin{cases} 
    \text{II-A-I$'$ : O.K.} \dotsb \; (**'') \\
    \text{II-A-II$'$} : Y':=Y+G \to 
     \begin{cases}
      \begin{gathered}
       \text{3, 4-a, } \\
       \text{5, 6-a} 
      \end{gathered} 
      \begin{cases}
       \text{II-A-I : O.K.} \dotsb \; (*'') \\
       \text{II-A-II} : Y'':=Y'+\phi' \to 
      \end{cases} \\
      \\
      \begin{gathered}
       \text{4-b} \\
       \text{6-b} 
      \end{gathered} : z'':=z'+\gamma'_1 w^{L'} \\
      \hspace*{0.7cm} \to 
      \begin{cases}
       \text{II-A-I$'$ : O.K.} \dotsb (**'') \\
       \text{II-A-II$'$} : Y'':=Y'+G' \to 
      \end{cases}
     \end{cases}
   \end{cases}\\
   \\
   \text{10-c} : z':=z+\gamma_1 w \to 
    \begin{cases}
     {z'}^3 Y^2 \in f_0 
      \begin{cases}
       \text{II-A-I$'$ : O.K.} \dotsb \; (**'') \\
       \text{II-A-II$'$} : Y':=Y+G \to 
        \begin{cases}
         \text{II-A-I : O.K.} \dotsb \; (*'') \\
         \text{II-A-II} : Y'':=Y'+\phi' \to 
        \end{cases}
      \end{cases}\\
     \text{otherwise} : \text{I-A} 
    \end{cases} 
  \end{cases}\\
  \\
  \text{II-B} 
  \begin{cases}
   \begin{gathered} 
    \text{1$'$, 2$'$, 3$'$, } \\
    \text{4-a$'$, 5$'$, 7$'$ }
   \end{gathered} 
   \begin{cases}
    \text{II-B-I : O.K.} \dotsb \; (*''') \\
    \text{II-B-II} : Y':=Y+\phi \; \to 
    \begin{cases}
     \text{II-B-I : O.K.} \dotsb \; (*''') \\
     \text{II-B-II} : Y'':=Y'+\phi' \; \to \; \dotsb 
    \end{cases}
   \end{cases} \\
   \\
   \text{4-b$'$} : w':=w+\gamma_1 z^K \; \to 
   \begin{cases} 
    \text{II-B-I$'$ : O.K.} \dotsb \; (**''') \\
    \text{II-B-II$'$} : Y':=Y+G \to 
     \begin{cases}
      \text{3$'$, 4-a$'$} 
      \begin{cases}
       \text{II-B-I : O.K.} \dotsb \; (*''') \\
       \text{II-B-II} : Y'':=Y'+\phi' \to 
      \end{cases} \\
      \\
      \text{4-b}' : w'':=w'+\gamma'_1 z^{K'} \\
      \hspace*{0.5cm} \to 
      \begin{cases}
       \text{II-B-I$'$ : O.K.} \dotsb (**''') \\
       \text{II-B-II$'$} : Y'':=Y'+G' \to 
      \end{cases}     
     \end{cases}
   \end{cases}
  \end{cases}
 \end{cases} $} \\
\\
   The procedures of combinations of ``a coordinate change $\longrightarrow$" 
must finish in finite times. Namely, taking suitable coordinate changes in 
finite times if necessary, the situation can be reduced to the case that 
$(\{f_0=0\} \cap H, P)$ is rational i.e. the case \\ $(*\dotsb*^{(m)})$ of (I) 
or (II), after all. Thus, there exist a local coordinate system and a weight 
such that $(f_0, P)$ is rational. Furthermore, the condition 
$ \frac12 + \beta^{(n)} + \gamma^{(n)} + \delta^{(n)} > 1 $ is satisfied under 
each coordinate system appearing in the above procedures. \\

   Thus it is enough to show Claim 3.5 for the proof of Theorem 2.7. \\
\\
{\it Proof of Claim 3.5}.\\
Step 1. 

  Since $(f_{\lambda},0)$ is an isolated singularity, it follows that: \\
\hspace*{2cm} $ z^k $ or $ y z^k $ or $ z^k w \in f_{\lambda}, $ \; and 
                \; $ w^l $ or $ y w^l $ or $ z w^l \in f_{\lambda}. $ \\
Taking suitable coordinate change \\
\hspace*{1cm} $ y' := y + b_1 z^{K_1} + c_1 w^{L_1}, \quad
                z' := z + c_2 w^{L_2}, \quad
                w' := w + b_2 z^{K_2} $ \\
for some sufficiently large $ K_i, \; L_j \in {\bold N} $ and some 
$ b_i, \; c_j \in {\bold C} $, we have $ {z'}^{k'}, \; {w'}^{l'} \in 
f_{\lambda} $ and $ \Gamma(f_{\lambda}) $ is the same as before except adding  
compact faces touching some coordinate planes.\\
\\
Step 2 ([A]-Thm. XXII, [AGV]-12.7, Kouchnirenko [K]-Thm. I).

  Let $ y^J, z^K, w^L \in F(y,z,w) \in {\bold C}[y,z,w] $ and $ V=\bigcup_{i=0}
^3 V_i $ be a decomposition of the three-dimensional region of the positive 
orthant below the Newton boundary $ \Gamma(F) \subset {{\bold R}_{\geq 0}}^3 
= \{ (Y,Z,W) \in {{\bold R}_{\geq 0}}^3 \} $ which satisfies: \\
\hspace*{1cm} $V_0$ is the three-dimensional simplicial cone which has vertex 
              $ (0,0,0), \; (j,0,0), $ \\ 
\hspace*{1cm} $ (0,k,0), \; (0,0,l) $ with $ j, k, l \in {\bold Q}_{>0} \; ; 
              \; j \leq J, \; k \leq K, \; l \leq L $,\\
\hspace*{1cm} $V_1$ has a vertex $ (J,0,0) $, \ 
              $V_2$ has a vertex $ (0,K,0) $, \ 
              $V_3$ has a vertex $ (0,0,L) $, \\
\hspace*{1cm} $ \dim_{\bold R}(V_i \cap V_j) \leq 2 $ for $ i \ne j $, \; and 
              $ S_{1 2} = S_{2 3} = S_{3 1} = \emptyset $ for \\
\hspace*{1cm} $ S_{i 1} := V_i \cap (YZ \text{-plane}), \; 
                S_{i 2} := V_i \cap (ZW \text{-plane}), \;
                S_{i 3} := V_i \cap (WY \text{-plane}) $.\\
(See {\sc Figure} 5.)
\begin{figure}[h]
\setlength{\unitlength}{1mm}
\begin{picture}(155,155)(-60,-60)
  \multiput(0,0)(1.5,0){50}{\circle*{0.2}}
  \multiput(0,0)(0,1.5){50}{\circle*{0.2}}
  \multiput(0,0)(-1,-1){40}{\circle*{0.2}}
  \put(75,0){\vector(1,0){10}}
  \put(0,75){\vector(0,1){10}}
  \put(-40,-40){\vector(-1,-1){10}}
  \put(-2,87){$Y$}
  \put(-55,-55){$Z$}
  \put(87,-2){$W$}
  \put(25,-10){\thicklines\line(1,3){10}}
  \put(35,20){\thicklines\line(3,-2){15}}
  \put(50,10){\thicklines\line(2,-1){10}}
  \put(60,5){\thicklines\line(3,-1){15}}
  \put(75,0){\thicklines\line(-5,-1){50}}
  \put(-10,25){\thicklines\line(3,2){15}}
  \put(5,35){\thicklines\line(1,0){15}}
  \put(20,35){\thicklines\line(-3,4){15}}
  \put(5,55){\thicklines\line(-1,4){5}}
  \put(0,75){\thicklines\line(-1,-5){10}}
  \put(-40,-40){\thicklines\line(2,1){50}}
  \put(10,-15){\thicklines\line(-2,1){10}}
  \put(0,-10){\thicklines\line(-1,1){10}}
  \put(-10,0){\thicklines\line(-1,2){5}}
  \put(-15,10){\thicklines\line(-1,-2){25}}
  \put(-15,10){\thicklines\line(1,3){5}}
  \put(20,35){\thicklines\line(1,-1){15}}
  \put(10,-15){\thicklines\line(3,1){15}}
  \multiput(10,-15)(-1.5,-0.5){25}{\circle*{0.2}}
  \multiput(-15,10)(-0.5,-1.5){25}{\circle*{0.2}}
  \multiput(25,-10)(1.5,0.5){20}{\circle*{0.2}}
  \multiput(55,0)(-1,1){20}{\circle*{0.2}}
  \multiput(20,35)(-1,1){20}{\circle*{0.2}}
  \multiput(0,55)(-0.5,-1.5){20}{\circle*{0.2}}
  \put(-5,75){$J$}
  \put(-40,-45){$K$}
  \put(75,-5){$L$}
  \put(-15,60){$j$}
  \put(-20,-40){$k$}
  \put(60,-15){$l$}
  \put(10,10){$V_0$}
  \put(10,55){$V_1$}
  \put(-40,-25){$V_2$}
  \put(55,10){$V_3$}
  \put(-10,60){\line(2,-1){9}}
  \put(60,-10){\line(-1,2){4.5}}
  \put(-20,-35){\line(-1,1){7}}
\end{picture}
\caption{ }
\end{figure} \\
Suppose that $(F,0)$ is an isolated singularity, then 
{\allowdisplaybreaks
\begin{align*}
\mu(F,0) & \geq 3! \left( \sum_{i=0}^3 V_i \right) 
           - 2! \left( \sum_{i=0}^3 \sum_{j=1}^3 S_{i j} \right)
           + 1!(J+K+L) - 1 \\
         & = jkl-(jk+kl+lj)+j+k+l-1 \\
         & \quad + \sum_{i=1}^3 \left( 6V_i-2\sum_{j=1}^3 S_{i j} \right) 
                 + (J-j) + (K-k) + (L-l) \\
         & \geq (j-1)(k-1)(l-1).
\end{align*}}
Step 3. 

  If $ 0 < \delta' < \delta \leq \gamma < \gamma' < 1 $ and 
  $ \gamma + \delta = \gamma' + \delta' $ 
  then 
  $$ \left( \frac{1}{\gamma}-1 \right) \left( \frac{1}{\delta}-1 \right) <
     \left( \frac{1}{\gamma'}-1 \right) \left( \frac{1}{\delta'}-1 \right). $$
In fact, let $ \gamma + \delta = \gamma' + \delta' =: 1/c, \; 
c \in {\bold R} $, \ then 
{\allowdisplaybreaks
\begin{align*}
 & \quad  \frac{1}{\gamma' \delta'}-\frac{1}{\gamma'}-\frac{1}{\delta'}-
   \left( \frac{1}{\gamma \delta}-\frac{1}{\gamma}-\frac{1}{\delta} \right) \\
 & = (c-1)(\gamma'-\gamma)(c(\gamma+\gamma')-1)/
           \gamma\gamma'(1-c\gamma)(1-c\gamma') > 0.
\end{align*}}
Step 4.

  From Step 1 - Step 3, if $ \gamma'+\delta'= 1/c \leq 1/6 $ then 
{\allowdisplaybreaks
\begin{align*}
   \mu(f_{\lambda},0) 
 & \geq (3-1) \left( \frac{1}{\gamma'}-1 \right) 
              \left( \frac{1}{\delta'}-1 \right) \\
 & \geq (3-1) \left( \frac{6}{c\gamma'}-1 \right) 
              \left( \frac{6}{c\delta'}-1 \right) \\
 & >    (3-1) \left( \frac{1}{\alpha_3}-1 \right)
              \left( \frac{1}{\alpha_4}-1 \right) = \mu(f,0),
\end{align*}}
a contradiction. This completes the proof of Claim 3.5 and Theorem 2.7.
\ \ \  Q.E.D.\\

\vspace{1cm}
\end{document}